\documentclass[twoside]{irmaems}
\usepackage{amsfonts}
\usepackage{amsmath}
\usepackage{epsf}
\usepackage{latexsym}
\usepackage[dvips]{graphicx}
\makeindex
\begin{document}
\newcommand{\aalpha}{{\overline{\alpha}}}
\renewcommand{\Im}{\mathop{\mathrm{Im}}\nolimits}
\newtheorem{definition}{Definition}
\newtheorem{proposition}{Proposition}
\newtheorem{conjecture}{Conjecture}
\newtheorem{corollary}{Corollary}
\title{Dual Teichm\"uller and lamination spaces.}
\author{V.V.Fock and A.B.Goncharov}
\address{}
\maketitle
\begin{abstract}
We survey explicit coordinate descriptions for two ($\mathcal A$ and $\mathcal X$) versions of Teichm\"uller and lamination spaces for open 2D surfaces, and extend them to the more general set-up of  surfaces with distinguished collections of points on the boundary. Main features, such as mapping class group action, Poisson and symplectic structures and others, are described in these terms. The lamination spaces are interpreted as the tropical limits of the Teichm\"uller ones. Canonical pairings between lamination and Teichm\"uller spaces are constructed.  
\end{abstract}

\section{Introduction.}
The Teichm\"uller space of a surface $S$ is the space of complex structures on $S$ modulo diffeomorphisms isotopic to the identity. For open surfaces this definition should be made more precise since we need to specify the behaviour of the complex structure at the boundary of the surface thus giving different versions of these spaces.

The space of measured laminations for a surface $S$ is defined by William Thurs\-ton \cite{Thurston} as the space of foliations of a closed subset of $S$ with a given transverse measure, modulo diffeomorphisms isotopic to the identity. A dense subset of the space of laminations is formed by foliations with finite number of leaves --- just collections of weighted closed curves without intersections and self intersections. For open surfaces there are different versions of such spaces. 

All versions of both Teichm\"uller and lamination spaces share the following properties: 

\begin{enumerate}
\item They are topologically trivial.
\item Isotopy classes of embeddings of surfaces induce canonical maps between these spaces.
\item There exist  canonical compactifications of the Teichm\"uller spaces by the quotients
 of the lamination spaces by multiplication by positive real numbers. 
\item There exists a canonical Poisson or degenerate symplectic structure (depending on a version) on these spaces.
\item A lamination with integral weights provides a function on the Teichm\"uller space. Such functions form a basis in the algebra of regular functions on Teichm\"uller space with respect to a certain algebraic structure on it.    
\end{enumerate}

The original approach to Teichm\"uller spaces used extensively functional analysis and was highly nonconstructive. The aim of the present paper is to give a description of all these spaces for open surfaces using elementary algebra and geometry and make all the properties obvious.   
  
We also consider generalisations of Riemann surfaces --- the so-called {\em ciliated surfaces} --- surfaces with a distinguished set of points on the boundary. Though at the beginning it requires a little more elaborate definitions, it provides us with more simple examples of Teichm\"uller spaces containing most features of the general cases. Namely, the configurations of points on the boundary of a disk or a cylinder fit perfectly into the picture.

We included a bit more technical Section \ref{Pairing} about the canonical pairing between lamination and Teichm\"uller spaces in the case of ciliated surfaces. 

The present exposition borrows a lot from \cite{F} and section 11 of \cite{FG1}. We tried to restrict ourselves to the purely geometric part of the picture and make it as elementary as possible. 
It can be used as an elementary introduction to more algebraic subjects, such as classical and quantum cluster varieties. To give the reader a 
further perspective, let us briefly mention some other sources and related topics. 

Explicit coordinate description of the Teichm\"uller spaces for an open surface goes back to William Thurston \cite{Thurston} and Robert Penner \cite{Penner}. This subject was developed further in \cite{F}.  
Quantum Teichm\"uller spaces were constructed in \cite{Chekhov-Fock} and independently by Rinat Kashaev in \cite{Kashaev}. Jorg Teschner \cite{Techner} proved that the quantum Teichm\"uller space acts on the space of conformal blocks of the Liouville conformal 
field theory, as was conjectured in \cite{F, Chekhov-Fock}. In \cite{FG1} the higher Teichm\"uller spaces were defined, and it was shown that they parametrise  certain discrete, faithful representations of the fundamental group of the surface to a split real simple Lie group $G$ of higher rank. These spaces  are closely related to the ones   studied by Nigel Hitchin \cite{Hitchin}. In particular, for $G=SL_3(\mathbb R)$ the corresponding higher Teichm\"uller space turns out \cite{FG3} to coincide with the space of real projective structures on $S$ studied by William Goldman and Suhyong Choi \cite{Goldman,Choi-Goldman}. In \cite{FG1, FG2} it was shown that the ${\cal A}$ and ${\cal X}$ versions of the Teichm\"uller and lamination spaces can be obtained as the positive real and tropical points of  certain {\em cluster ${\cal A}$- and ${\cal X}$-varieties}.  The closely related objects,  cluster algebras, were introduced by Sergey Fomin  and Andrei Zelevinsky \cite{FZ1}, and studied by them and Arkady Berenstein, Mikhail Gekhtman, Mikhail Shapiro, Alek Vainshtein (\cite{Zelevinsky} and references therein). As a result of these developments, the cluster theory was enriched by new examples as well as new features (such as duality, Poisson structure and quantisation, relations to the algebraic $K$-theory and the dilogarithm).  

In particular the canonical pairings from section 6, in the special case of a ciliated disc, can be viewed as the canonical pairings for cluster ${\cal A}$ and ${\cal X}$-varieties of finite type  $A_n$, predicted by the general duality conjectures \cite{FG2}.                                                                                  

We are very grateful to Athanase Papadopoulos and Vladimir Turaev, who suggested to write this paper for the Handbook on Teichm\"uller spaces. The first author is grateful to the Fields Institute and Brown University, where this paper was mainly written, for hospitality. The authors were supported by the NSF grants DMS-0400449 (the second author) and CRDF 2622;2660 (the first author). We are very indebted also to Andrei Levin, Yuri Neretin, Stepan Orevkov,  Athanase Papadopoulos, Robert Penner and especially to Aleksei Rosly for very valuable discussions and for Guillaume Th\'eret for carefully reading the text and correcting many typos and mistakes. 

\section{Surfaces and triangulations.}
In this section we shall briefly recall the basic facts about triangulations of 2D surfaces.

A {\em ciliated surface}\index{surface!ciliated} is a 2D compact oriented surface with boundary and with a finite set of marked points on the boundary called {\em cilia}\index{cilium}. 

A boundary component without cilia is called a {\em hole}.

A {\em triangulation}\index{triangulation} $\Gamma$ of a ciliated surface is a decomposition of the surface with contracted holes into triangles such that every vertex of a triangle is either a cilium or a shrunk hole.

\begin{figure}[h]
\center
\hspace{8mm}\begin{minipage}{3.2cm}
\includegraphics[scale=0.3]{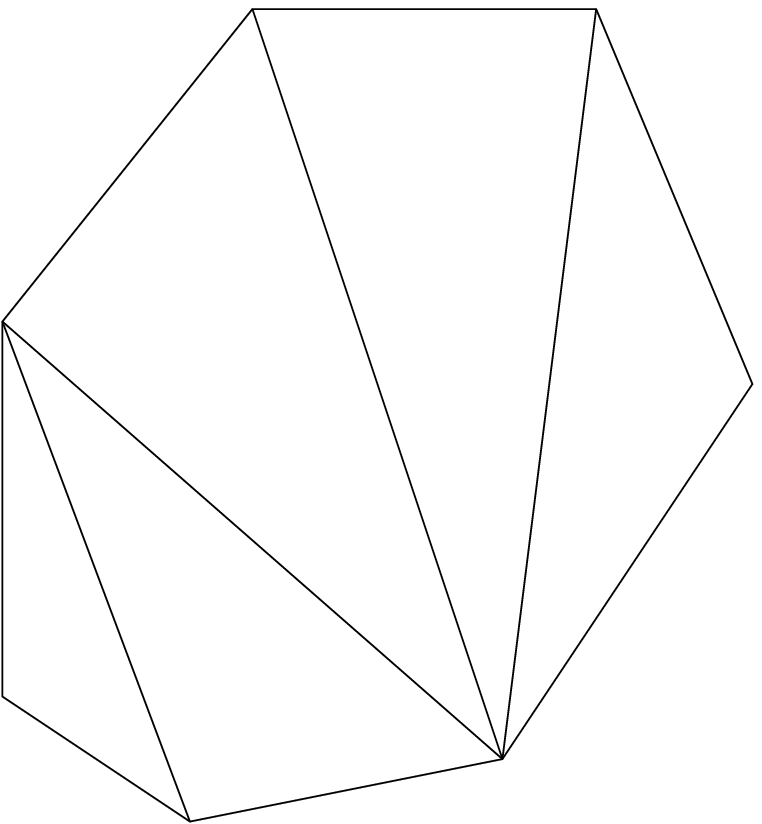}\par
\noindent {\textbf 1.}$g=0,P=(7)$
\end{minipage}
\begin{minipage}{3.2cm}
\includegraphics[scale=0.3]{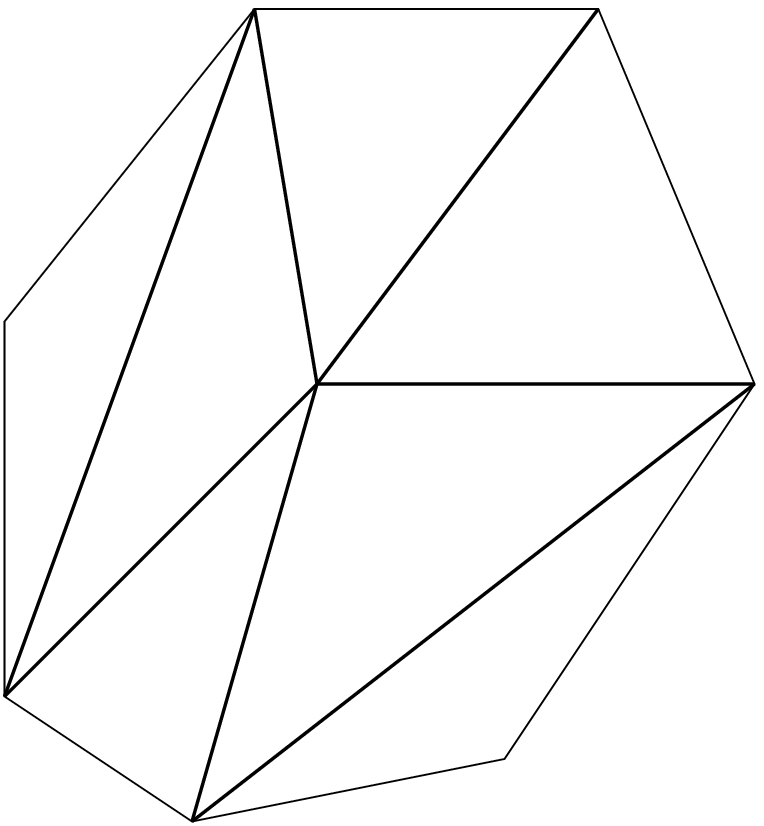}\par
\noindent {\textbf 2.}$g=0,P=(7,0)$
\end{minipage}
\begin{minipage}{3.2cm}
\includegraphics[scale=0.3]{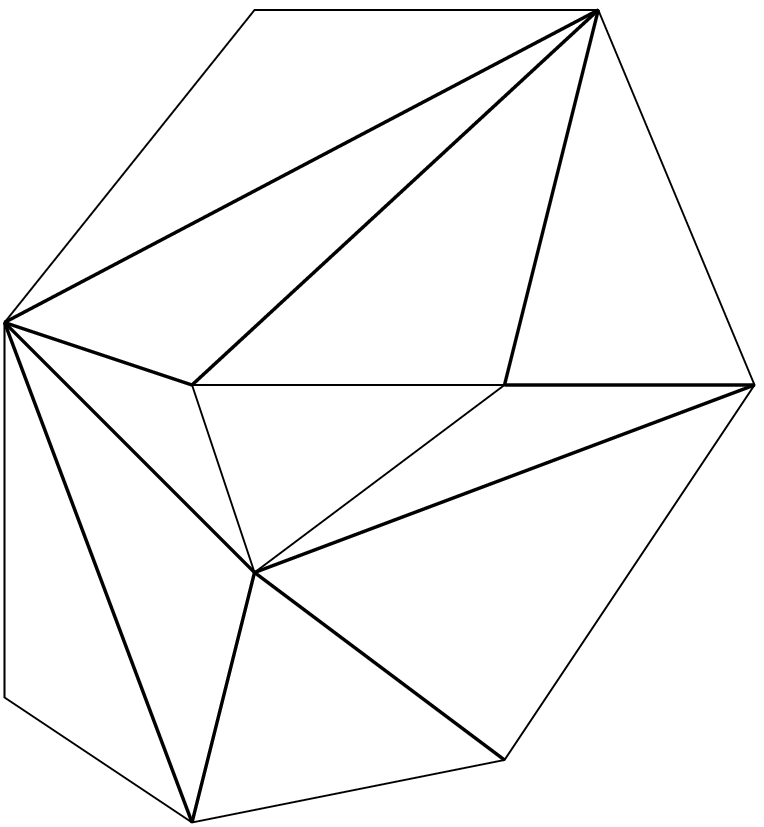}\par
\noindent {\textbf 3.}$g=0,P=(7,3)$
\end{minipage}
\vspace{3mm}

\begin{minipage}{6cm}
\includegraphics[scale=0.2]{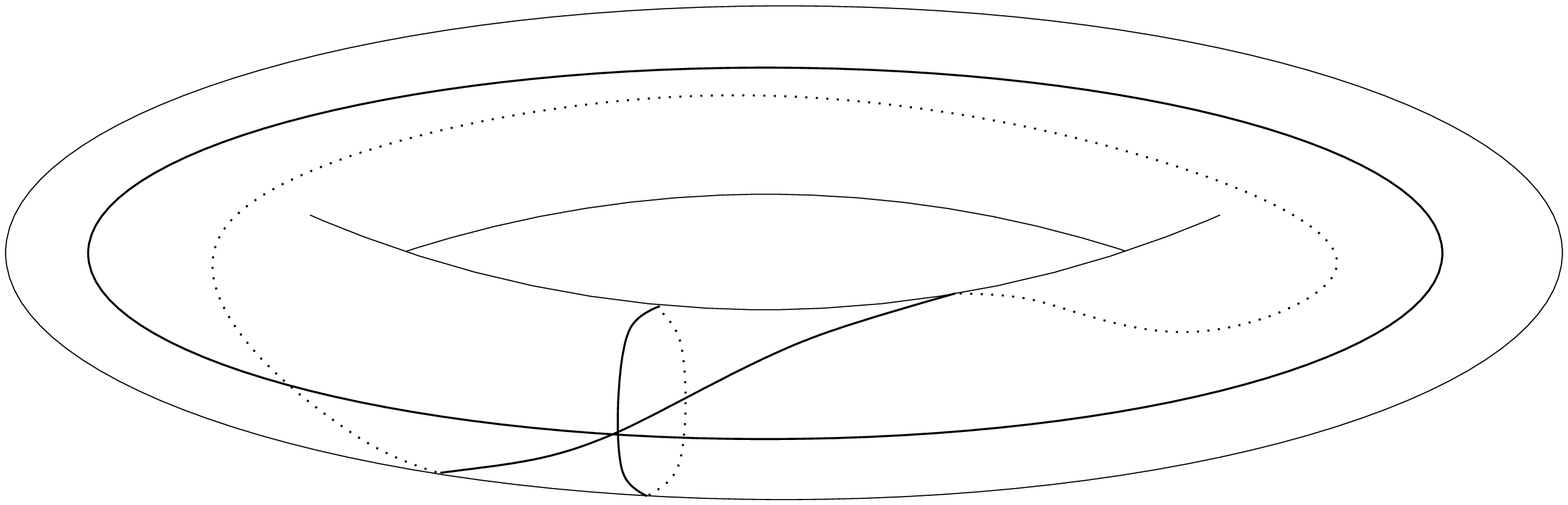}\par\vspace{1mm}
\hspace{1.5cm}{\textbf 4.}$g=1,P=(0)$
\end{minipage}
\caption{Examples of triangulations of ciliated surfaces. Internal edges are thick.}
\label{example}
\end{figure}

The edges of the triangulation belonging to two (resp. one) triangle are called {\em internal}\index{edge!internal} (resp. {\em external}\index{edge!external}). The external edges form the boundary of the surface. Denote by $F(\Gamma),E(\Gamma),E_0(\Gamma),V(\Gamma)$ the set of triangles, edges, external edges and vertices of the triangulation, respectively.

The vertices of the triangulation with the link homeomorphic to a circle (resp. a segment) correspond to holes (resp. cilia). Let us call the former type of vertices the {\em holes}\index{hole} as well.

Topologically a ciliated surface is defined by its genus $g$ and a finite collection $P=(p_1,\ldots,p_s)$, where $s$ is the number of boundary components and $p_i$ is the number of cilia on the $i$-th component. There is no canonical order on the set of boundary components. Denote the number of holes by $h$, the total number of cilia by $c=\sum_i p_i$, and the number of internal edges by $n$.

Observe, that the numbers of faces, vertices and edges, both external and internal, are determined by the topology.
\begin{enumerate}
\item  $\sharp V(\Gamma)=h+c$
\item  $\sharp E_0(\Gamma)=c$
\item  $\sharp E(\Gamma)=6g-6+3s+2c$
\item  $n = 6g-6+3s+c$
\item  $\sharp F(\Gamma)=4g-4+2s+c$
\end{enumerate}
The the first two statements are obvious. The others follow from the expression for the Euler characteristics of the surface:
$$\sharp F(\Gamma)-\sharp E(\Gamma)+ \sharp V(\Gamma)=2-2g+h-s$$
and the expression
$$3\sharp F(\Gamma)= 2\sharp E(\Gamma)-\sharp E_0(\Gamma)$$
of the fact that each triangle has three sides.

The topology of the triangulation can be encoded into a skew-symmetric matrix $\varepsilon^{\alpha\beta}$, where $\alpha,\beta \in E(\Gamma)$, defined as
$$
\varepsilon^{\alpha\beta} = \sum_{i\in F(\Gamma)}<\alpha, i, \beta>, 
$$
where $<\alpha, i, \beta>$ is equal to $+1$ (resp. $-1$) if both $\alpha$ and $\beta$ are sides of the triangle $i$ and $\alpha$ is in the counterclockwise (resp. clockwise) direction from $\beta$ with respect to their common vertex. Otherwise $<\alpha, i, \beta>$ is equal to zero. The entries of the matrix $\varepsilon^{\alpha\beta}$ might have five possible values: $0,\pm 1,\pm 2$.

\noindent
\begin{minipage}{2.5cm}
\unitlength 0.8pt
\begin{picture}(0,97)(-10,40)
\put(10,60){\line(1, 1){20}} 
\put(10,60){\line(1,-1){20}} 
\put(10,60){\line(1,0){40}} 
\put(50,60){\line(-1,1){20}}
\put(50,60){\line(-1,-1){20}}
\put(10,118){\line(1,1){20}}
\put(10,118){\line(1,-1){20}}
\put(50,118){\line(-1,1){20}}
\put(50,118){\line(-1,-1){20}}
\put(30,138){\line(0,-1){40}}
\put(30,89){\makebox(0,0)[cc]{$\updownarrow$}}
\end{picture}\par
\refstepcounter{figure}\label{flip}
\noindent Figure \ref{flip}. A flip.
\end{minipage}
\hfill
\begin{minipage}{9cm}
The number of triangulations of a given ciliated surface is infinite except for $g=0$ and for either $P=(k)$ or $P=(k,0)$. However one triangulation can be obtained from another one by a sequence of elementary moves called {\em flips}\index{flip} or {\em Whitehead moves}. One triangulation is a flip of another if it is obtained by removing one internal edge and replacing it by another diagonal of the arising quadrilateral (fig. \ref{flip}). A flip can be done with any internal edge unless it belongs to one triangle of the triangulation.
\end{minipage}\vspace{0.5mm}
Observe that there is a canonical correspondence between the edges of a triangulation and the edges of a flipped one. We shall use the notation $\alpha'$ for the edge corresponding to the edge $\alpha$. However one should be aware that compositions of a sequence of flips may restore the original triangulations but a nontrivial correspondence between the edges. An example of this phenomenon can be illustrated by fig. \ref{pentagon} where the composition of five flips does not change the triangulation of the pentagon, but the canonical correspondence interchanges the diagonals. 

It is a simple exercise to check that the matrix $\varepsilon^{\alpha'\beta'}$ encoding the combinatorics of a triangulation after a flip in the edge $\gamma$ is given by the formula:

\begin{equation}\label{epsilon-flip}
\varepsilon^{\alpha'\beta'}=
\left\{\begin{array}{ll}
     -\varepsilon^{\alpha\beta}&\mbox{if }\alpha=\gamma \mbox{ or } \beta=\gamma\\
     \varepsilon^{\alpha\beta} + \frac{1}{2}(\varepsilon^{\alpha\gamma}|\varepsilon^{\gamma\beta}|+|\varepsilon^{\alpha\gamma}|\varepsilon^{\gamma\beta})&
\mbox{otherwise}
\end{array}  
 \right.
\end{equation}

To give some definitions we shall need to fix the orientation of the holes of the surface. By that we mean fixing orientation of the boundary components in a way not necessarily induced by the orientation of the surface. We say that the orientation of a hole is {\em positive} (resp. {\em negative}) if it agrees (resp. disagrees) with the orientation of $S$. 

The {\em mapping class group}\index{mapping class group} $\mathcal D(S)$ of a ciliated surface $S$ is the group of connected components of the diffeomorphisms of $S$ preserving the set of cilia. In our examples of fig. \ref{example} the mappings class groups are, respectively $\mathbb Z/7\mathbb Z$, $\mathbb Z/7\mathbb Z$, $\mathbb Z/3 \mathbb Z\times \mathbb Z/7 \mathbb Z$ and $PSL(2,\mathbb Z)$. Observe, that the mapping class group is finite if and only if the number of triangulations is finite.

 In the appendix \ref{combinatorial} we give a combinatorial description of the mapping class group in terms of the surface triangulations.

\section{Laminations.}\index{lamination}

 Taking into account that the reader may be unfamiliar with the Thurston's notion of a measured lamination \cite{Thurston}, we are going to give all definitions here in the form which is almost equivalent to the original one (the only difference is in the treatment of the holes, punctures and cilia), but more convenient for us. The construction of coordinates on the space of laminations we are going to describe is a slight modification of Thurston's "train tracks" (\cite{Thurston}, section 9).

It seems worth mentioning here that the definitions of measured laminations are very similar to the definitions of the singular homology groups, and is in a sense an unoriented version of the latter ones.

There are two different ways to define the notion of measured laminations for surfaces with boundary, which are analogous to the definition of homology group with compact and closed support, respectively.

Measured laminations are certain collection of weighted curves on the surface. By a {\em curve} here we mean a curve  without self-intersections either closed or connecting two points on the boundary, disjoint from cilia and considered up to homotopy within the class of such curves. We call a curve {\em special} if it is retractable to a hole or to a interval on the boundary containing exactly one cilium. We call a curve {\em contractible} if it can be retracted to a point within this class of curves. In particular a nonclosed curve is contractible if and only if  it is retractable to a segment of the boundary containing no cilia.

\subsection{Unbounded measured laminations.}\index{lamination!unbounded}

\paragraph{Definition}
{\em A rational unbounded measured lamination} or an {\em $\mathcal X$-lamination} on a 2-dimensional ciliated surface is given by orientations of some boundary components and a homotopy class of a collection of finite number of non-self\-in\-ter\-sec\-ting and pairwise nonintersecting curves with positive rational weights either closed or ending at the boundary and disjoint from cilia, subject to the following equivalence relations:
\begin{enumerate}
\item A lamination containing a contractible or a special curve is equivalent to the lamination with this curve removed.
\item A lamination containing two homotopy equivalent curves of weights $u$ and $v$, respectively, is equivalent to the lamination with one of these curves removed and with the weight $u+v$ on the other.
\item Oriented boundary components are holes unless they are disjoint
  from the curves.
\end{enumerate}

 The set of all rational unbounded laminations on a given surface $S$ is denoted by $\mathsf T^x(S,\mathbb Q)$. This space has a natural subset, given by collections of curves with integral weights. This space is denoted by $\mathsf T^x(S,\mathbb Z)$. We shall also omit the arguments if they are clear from the context.

There is an  action of the multiplicative group of positive numbers on unbounded laminations $\mathbb Q_{> 0}\times \mathsf T^x \rightarrow \mathsf T^x $ given by multiplication of the weights of all curves by a fixed positive number. 

For every hole $\rho$ we associate a number $\mathsf r^\rho$, which is the total weight (resp. minus the total weight) of curves entering the hole  if the orientation of the hole is positive (resp. negative).

\begin{figure}[h]
\center
\includegraphics[scale=0.15]{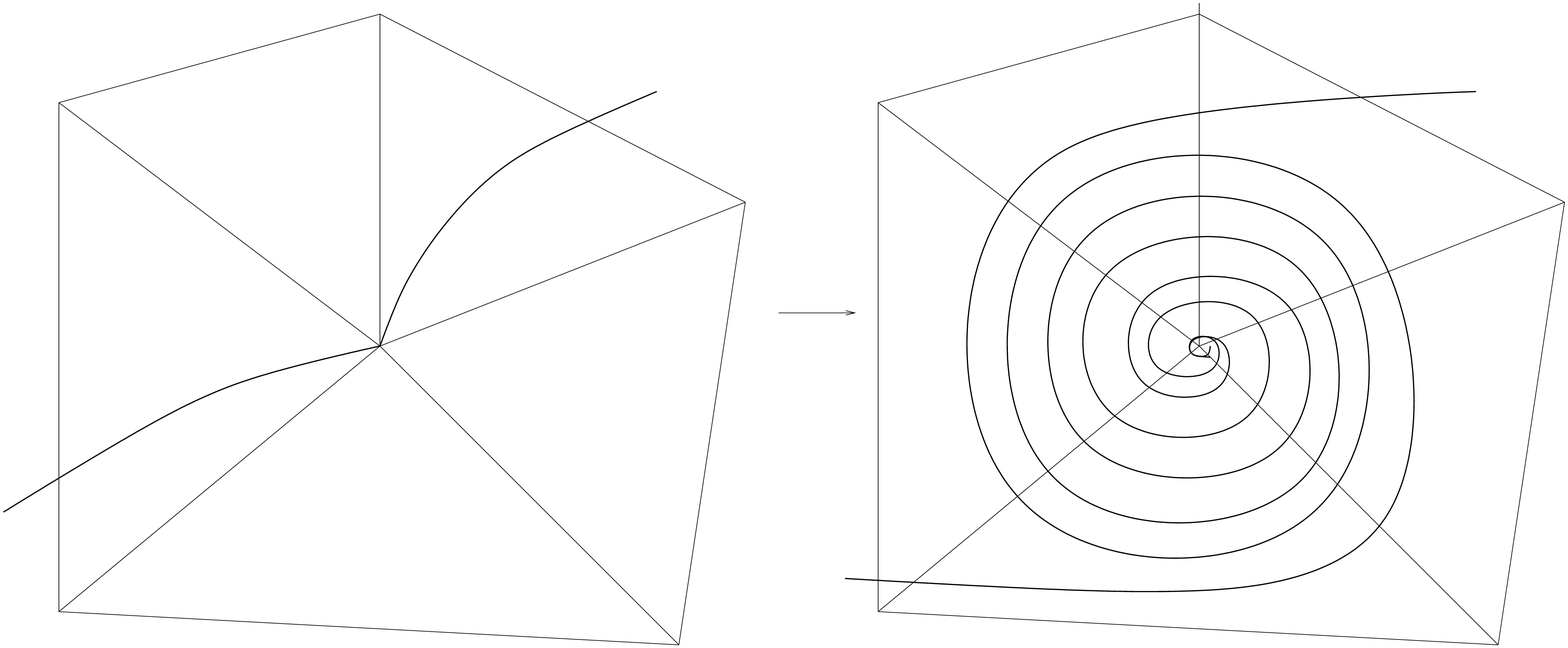}
\caption{Twisting curves incident to a boundary component without cilia.}
\label{twist}
\end{figure}

\paragraph{Construction of coordinates.}

Let us first rotate each boundary component without cilia infinitely many times in the direction prescribed by the orientation as shown on fig. \ref{twist}.

\begin{minipage}{5cm}
\begin{picture}(0,140)(0,-20)
\includegraphics[scale=0.3]{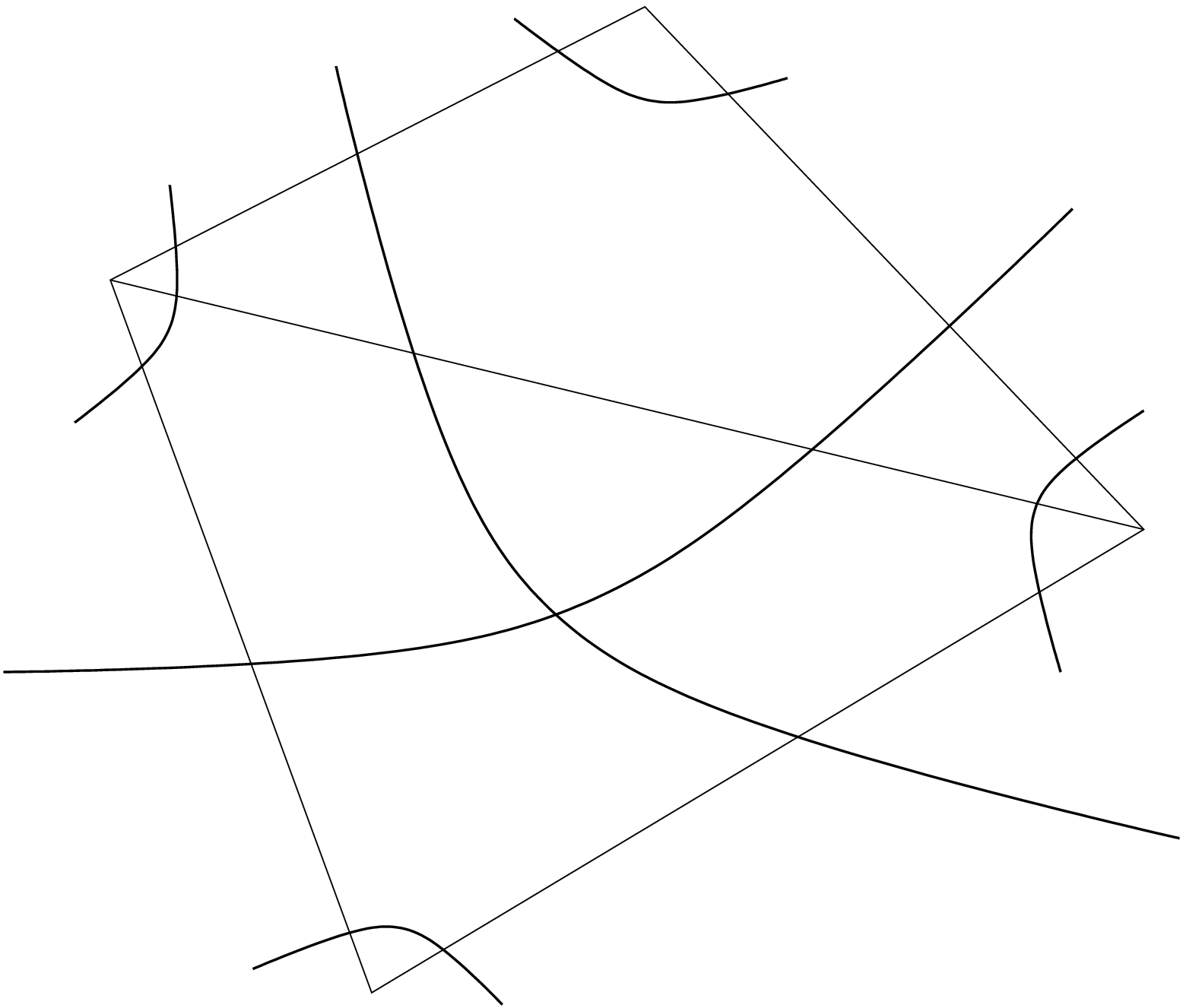}
\put(3,56){\makebox(0,0)[cc]{$C$}}
\put(-63,123){\makebox(0,0)[cc]{$B$}}
\put(-93,-6){\makebox(0,0)[cc]{$D$}}
\put(-130,83){\makebox(0,0)[cc]{$A$}}
\put(-120,43){\makebox(0,0)[cc]{negative}}
\put(-13,30){\makebox(0,0)[cc]{positive}}
\put(-63,73){\makebox(0,0)[cc]{$\alpha$}}
\end{picture}
\refstepcounter{figure}\label{intersection}
\noindent Figure \ref{intersection}. Positive and negative intersection points.
\end{minipage}
\hfill
\begin{minipage}{6.2cm}
  Then deform the curves in such a way that they do not cross any edge of the triangulation consecutively in the opposite directions.  Consider the quadrilateral formed by two triangles sharing an edge $\alpha$. Denote its vertices by $ABCD$ in the clockwise direction with respect to the orientation, starting from a vertex of $\alpha$. We call an intersection point of the edge $\alpha$ and a curve of the lamination {\em positive} if it belongs to the segment of the curve connecting $AB$ and $CD$ and {\em negative} if it belongs to the segment of the curve connecting $BC$ and $AD$. Since the curves of the lamination do not intersect, having both negative and positive intersection points on the same edge is impossible. Notice that we do not assign a sign, say, to a curve intersecting $BC$ and $CD$. \vspace{3pt}
\end{minipage}

Now assign to each edge $\alpha$  the sum over the positive intersection points of weights of the respective curves minus the sum over the negative intersection points of weights of the respective curves. The collection of these numbers, one for each internal edge of $\Gamma$, is the desired set of coordinates.

Observe that although the curves spiralling around the holes without cilia intersect some edges infinitely many times, only finite number of these intersection points are positive or negative, and therefore all numbers on edges obtained is this way are finite.

Now we need to prove that these numbers are coordinates indeed. We shall do it by describing an inverse construction. 

Note that if we are able to construct a lamination corresponding to the set of numbers $\{\mathsf x^\alpha\}$, we can also do it for the set $\{u \mathsf x^\alpha\}$ for any rational $u \ge 0$.  Therefore we can reduce our task to the case when all numbers on edges are integral.

\begin{figure}
\center
\includegraphics[scale=0.25]{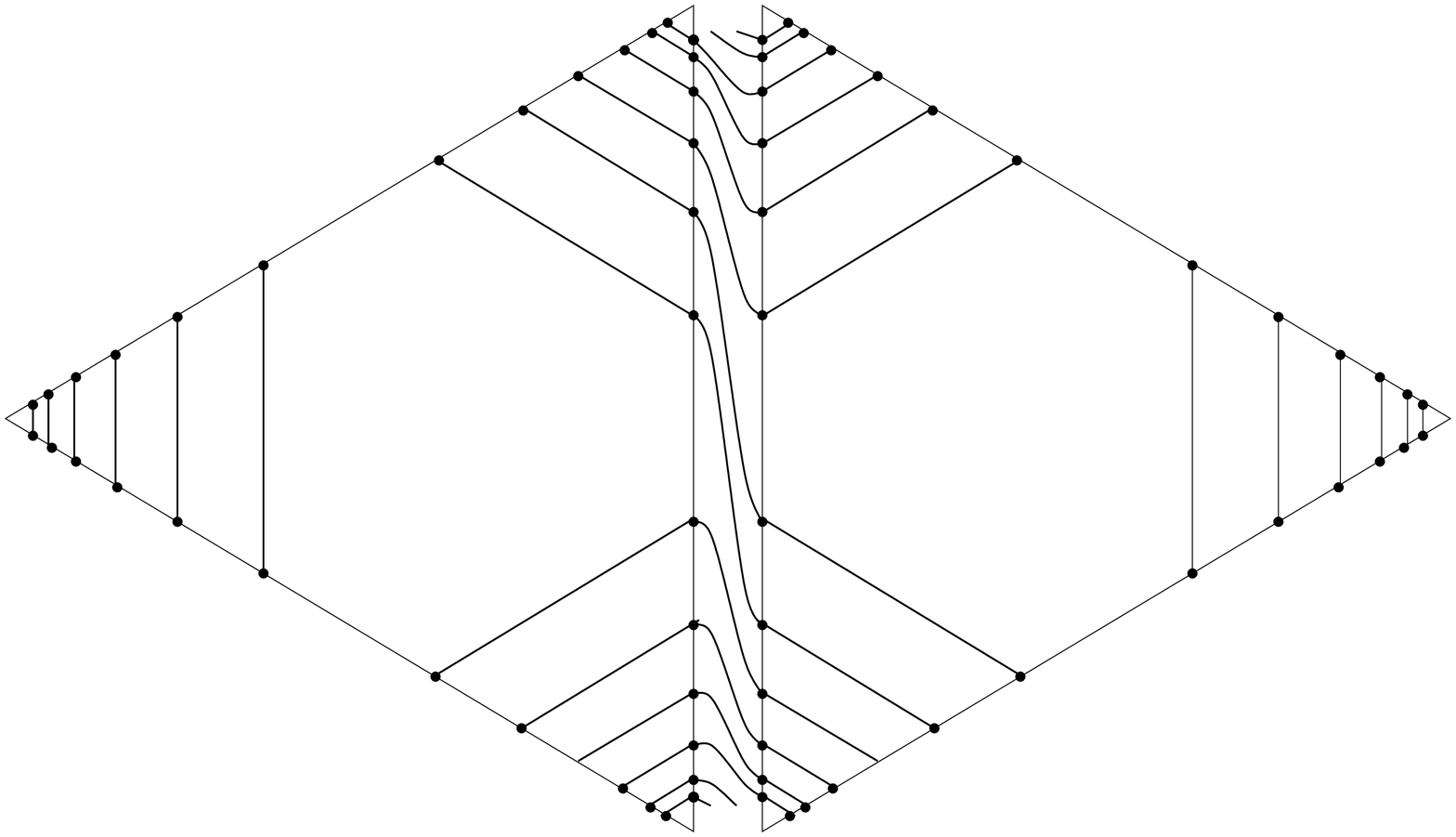}
\caption{Gluing triangles for $\mathsf x^\alpha=2$.}
\label{reconstructionX}
\end{figure}

Take a triangle, parametrise each side of it by $\mathbb R$ respecting orientation induced by the orientation of the triangle, and connect points with parameter $i\in \frac12 + \mathbb Z_{\geq 0}$ on one side to the point with parameter $-i$ on the next side in the clockwise direction. It can be done by curves without mutual intersections. Now replace every triangle of our triangulation with such triangles with curves in such a way that the point with parameter $i$ on one side of an edge $\alpha$ is glued to the point with parameter $\mathsf x^\alpha-i$ on the other side (fig. \ref{reconstructionX}). 

 Observe that although we have started with infinite bunches of curves, the resulting lamination is finite: all these curves glue together into a finite number of connected components and possibly infinite number of special curves. Indeed, any connected component either intersects positively or negatively at least one edge or is closed. Since the total number of positive or negative intersection points  $\sum_{\alpha \in E(\Gamma)- E_0(\Gamma)}\vert\mathsf x^\alpha\vert$ is finite, the resulting lamination contains no more than this number of connected components. In particular the number of connected components equals this number provided all numbers $\mathsf x_\alpha$ are all nonpositive or all nonnegative.

\paragraph{Properties of the coordinates.}
 The constructed coordinates on the space of laminations correspond to a particular choice of the triangulation.  If we change the triangulation, the corresponding coordinates change. The rule how the coordinates change under a flip of the triangulation is given by an explicit formula:
\begin{equation}\label{xl-flip}
\mathsf x'^{\beta'}=\left\{\begin{array}{ll}
-\mathsf x^\alpha&\mbox{ if }\beta=\alpha\\
\mathsf x^\beta+\varepsilon^{\alpha\beta}\max(0,\mathsf x^\alpha) &\mbox{ if }\varepsilon^{\alpha\beta}\geq 0\\
\mathsf x^\beta+\varepsilon^{\alpha\beta}\max(0,-\mathsf x^\alpha) &\mbox{ if }\varepsilon^{\alpha\beta}\leq 0
\end{array}
\right.
\end{equation}

\begin{figure}[h]
\unitlength 1mm
\begin{picture}(100,50)(4,95)
\put(70,120){\line(1, 1){20}} \put(76,108){\makebox(0,0)[cc]{$\mathsf x^4-\max(-\mathsf x^0,0)$}}
\put(70,120){\line(1,-1){20}} \put(76,132){\makebox(0,0)[cc]{$\mathsf x^1+\max(\mathsf x^0,0)$}}
\put(70,120){\line(1,0){40}} \put(90,120) {\makebox(0,0)[bc]{$-\mathsf x^0$}}
\put(110,120){\line(-1,1){20}}\put(105,108){\makebox(0,0)[cc]{$\mathsf x^3+\max(\mathsf x^0,0)$}}
\put(110,120){\line(-1,-1){20}}\put(105,132){\makebox(0,0)[cc]{$\mathsf x^2-\max(-\mathsf x^0,0)$}}
\put(10,120){\line(1,1){20}} \put(18,108){\makebox(0,0)[cc]{$\mathsf x^4$}}
\put(10,120){\line(1,-1){20}} \put(18,132){\makebox(0,0)[cc]{$\mathsf x^1$}}
\put(50,120){\line(-1,1){20}}\put(42,108){\makebox(0,0)[cc]{$\mathsf x^3$}}
\put(50,120){\line(-1,-1){20}}\put(42,132){\makebox(0,0)[cc]{$\mathsf x^2$}}
\put(30,140){\line(0,-1){40}}\put(32,120) {\makebox(0,0)[cc]{$\mathsf x^0$}}
\put(60,120){\makebox(0,0)[cc]{$\longleftrightarrow$}}
\end{picture}
\caption{Changing coordinates under a flip of the triangulation. Only the changing coordinates are shown here, the numbers on the other edges remain unchanged.}
\label{xflipl}
\end{figure}
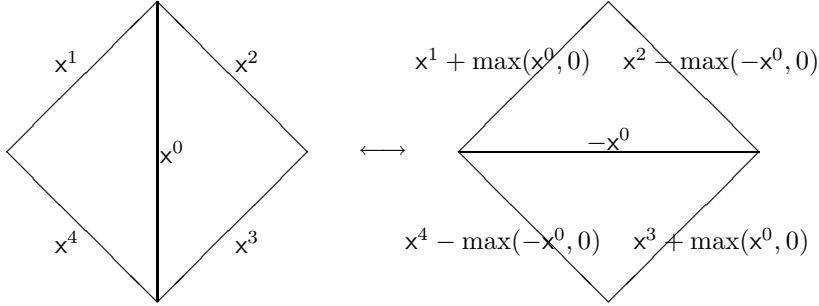

If all edges of the quadrilateral taking part in the flip are different, the change of coordinates can be shown by the graphical rule fig.  \ref{xflipl}.

An ${\mathcal X}$-lamination is {\em integral}\index{lamination!integral} if all its coordinates are integral. 

If $\rho$ is a hole the canonical map $\mathsf r^\rho\colon \mathsf T^x \rightarrow \mathbb Q$ is given by $\mathsf r^\rho= \sum_\alpha \mathsf x^\alpha$, where the sum is taken over all edges incident to $\rho$ (the edge $\alpha$ is taken twice if both ends of $\alpha$ are incident to $\rho$).

\subsection{Bounded measured laminations.}\index{lamination!bounded}

\paragraph{Definition.}{\em A rational bounded measured lamination} or a {\em rational $\mathcal A$-la\-mi\-na\-tion} on a 2-dimensional ciliated surface is a homotopy class of a collection of finite number of self- and mutually nonintersecting unoriented curves  either closed or connecting two points of the boundary disjoint from cilia with rational weights and subject to the following conditions and equivalence relations. 

\begin{enumerate}

\item Weights of all curves are positive, unless a curve is special.

\item A lamination containing a curve of weight zero is considered to be
equivalent to the lamination with this curve removed.

\item A lamination containing a contractible curve is considered to be equivalent to the lamination with this curve removed.

\item A lamination containing two homotopy equivalent curves with weights $u$
and $v$ is equivalent to the lamination with one of these curves
removed and with the weight $u+v$ on the other.

\end{enumerate}

Recall that a contractible curve in this context is either closed or retractable to a segment of the boundary without cilia.

 The set of all rational bounded laminations on a given surface $S$ is denoted by $\mathsf T^a(S,\mathbb Q).$ This space has a natural subset called the set of {\em integral}\index{lamination!integral} bounded laminations, consisting of laminations with integral weights.  It is denoted by $\mathsf T^a(S,\mathbb Z)$. We shall omit the arguments of $\mathsf T^a$ if they are clear from the context or if the statement is valid for any value of the arguments.

The space $\mathsf T^a(S)$ has a subspace $\mathsf T^a_0(S)$ consisting of the laminations such that for any segment of the boundary between two cilia the total weight of curves ending at it vanishes.

There is a canonical action of the multiplicative group of positive numbers on bounded lamination $\mathbb Q_{> 0}\times \mathsf T^a \rightarrow \mathsf T^a $ given by multiplication of weights of all curves by a fixed positive number.

For every hole $\rho$ there is an action $\mathsf r_\rho$ of the additive group of rational numbers $\mathsf r_\rho\colon\mathbb Q\times\mathsf T^a\rightarrow \mathsf T^a$ by adding a closed loop with certain weight around the hole. Similarly, for every cilium $\rho$ there is an action (also denoted by $\mathsf r_\rho$) of the additive group $\mathbb Q$  by adding a special curve about this cilium with a certain weight. 

For every hole $\rho$ there is a map $\mathsf A_\rho \colon \mathsf T^a\rightarrow \mathbb Q$ (called a {\em collar map}\index{collar map} \cite{Papadop}, where it was defined and studied) given by the total weight of the loops surrounding the hole. 

There is a canonical map $\mathsf p\colon \mathsf T^a(S)\rightarrow \mathsf T^x(S)$ since every bounded lamination can be considered as an unbounded one. This map obviously commutes with the multiplication by positive numbers and sends integral laminations to the integral ones. 

\paragraph{Construction of coordinates.}

Suppose we are given a triangulation $\Gamma$ of a ciliated surface $S$.  We are going to assign, for a given lamination, rational numbers on edges of $\Gamma$ and show that these numbers are global coordinates on the space of laminations.

Deform the curves of the lamination in such a way that every curve intersects every edge at the minimal possible number of points. Assign to each edge $\alpha$  a half of the sum over the intersection points of weights of the respective curves. The collection of these numbers, one for each edge of $\Gamma$, is the desired set of coordinates. For the space $\mathsf T_0^a(S)$ the numbers assigned in this way to any external edge vanishes, thus only numbers assigned to internal edges serve as coordinates for this space.

\paragraph{Reconstruction.}
Now we need to prove that these numbers  are coordinates indeed. For this purpose we just describe an inverse construction, which gives a lamination, starting from numbers on edges.

 Let us first fix a pair of rationals $u$ and $v$ and introduce the numbers $\{\tilde {\mathsf a}_\alpha|\alpha \in E(\Gamma)\}$ by 
$$\tilde {\mathsf a}_\alpha= u\mathsf a_\alpha+v
$$ 
By choosing $u$ and $v$ one can make numbers $\{\tilde {\mathsf a}_\alpha\}$ to be positive integers and to satisfy the triangle inequality for every triangle of the triangulation with edges $\alpha,\beta,\gamma$:
\begin{equation}
\vert \tilde {\mathsf a}_\alpha-\tilde {\mathsf a}_\beta \vert \le \tilde {\mathsf a}_\gamma \le  \tilde {\mathsf a}_\alpha+\tilde {\mathsf a}_\beta.
\end{equation}

\noindent
\begin{minipage}[c]{5.5cm}
~~\includegraphics[scale=0.25]{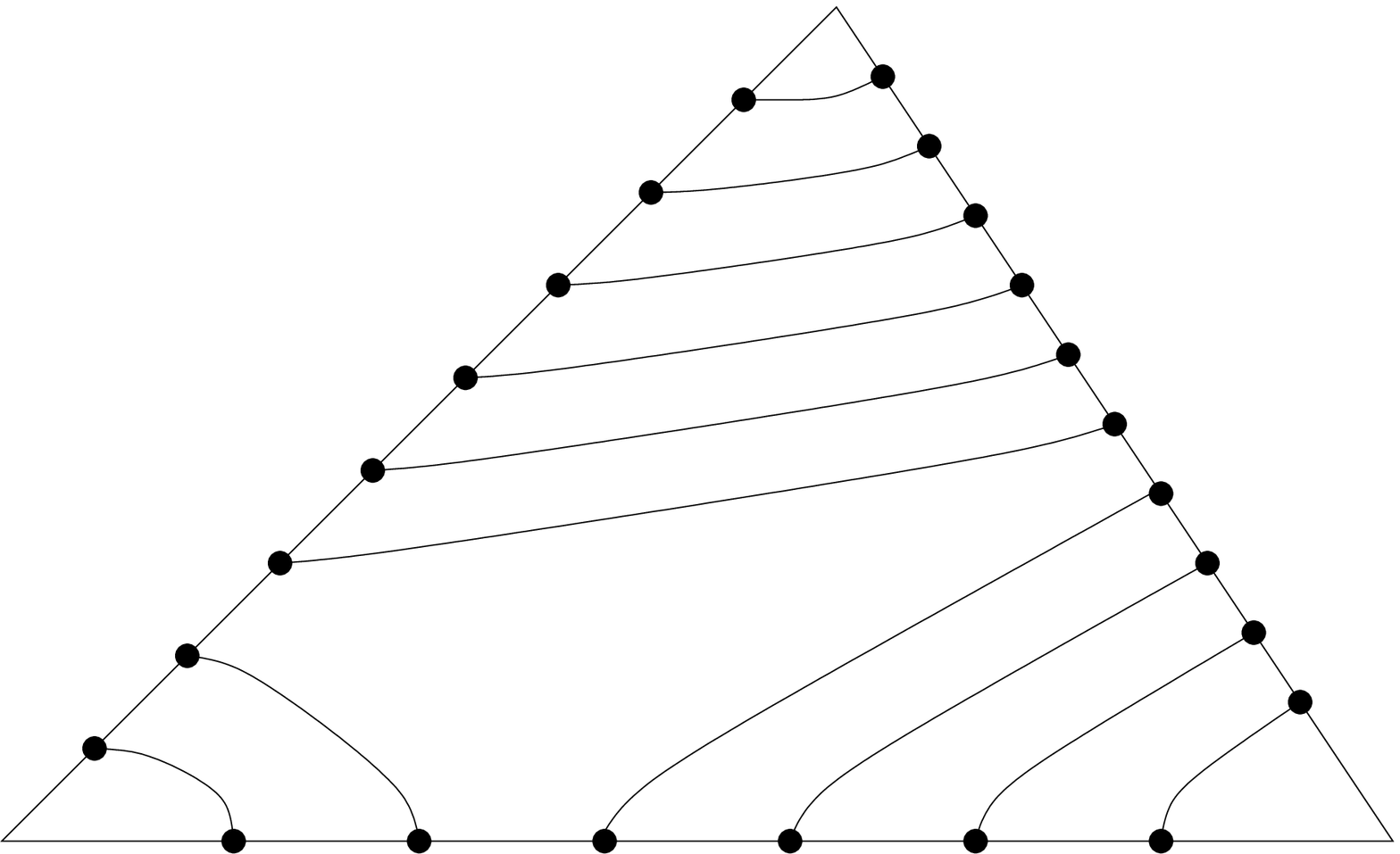}\\
\refstepcounter{figure}\label{vertex}
\noindent Figure \ref{vertex}. $(\tilde {\mathsf a}_\alpha,\tilde {\mathsf a}_\beta,\tilde {\mathsf a}_\gamma)=(3,4,5)$
\end{minipage}\hfill
\begin{minipage}{6.2cm}Now mark $2\tilde {\mathsf a}_\alpha$ points on each edge $\alpha$ and join the points on the sides of each triangle pairwise by nonintersecting and non-selfintersecting lines connecting points on different sides an passing inside the triangle. One can easily see that once the triangle inequality is satisfied this can be done in a unique way up to homotopy (fig.\ref{vertex}). Now take each hole and each cilium and add a closed curve surrounding them with weight $-v$. Then divide the weights of all curves by $u$. The construction is finished. The proof that the construction is indeed inverse to the previous one is self evident.\vspace{3pt}
\end{minipage}

\paragraph{Properties of the coordinates.}
 Just as for the case of $\mathcal X$-laminations, the constructed coordinates on the space of laminations correspond to a particular choice of the triangulation.  The rule how the coordinates change under a flip of the triangulation is given by the formula:
\begin{equation}\label{al-flip}
\mathsf a'_{\beta'}=\left\{\begin{array}{lcl}
\max(\sum\limits_{\delta|\varepsilon^{\beta\delta}>0}\varepsilon^{\beta\delta}\mathsf a_\delta,-\sum\limits_{\delta|\varepsilon^{\beta\delta}<0}\varepsilon^{\beta\delta}\mathsf a_\delta)-{\mathsf a_\alpha}&\mbox{ if }\beta=\alpha\\
\mathsf a_\beta &\mbox{ if }\beta\neq\alpha\\
\end{array}
\right.
\end{equation}

  Since any triangulation change is a composition of flips, these rules allow to express the coordinate change for any triangulation change.

\begin{figure}[h]
\unitlength 1mm
\begin{picture}(100,50)(4,95)
\put(70,120){\line(1, 1){20}} \put(76,108){\makebox(0,0)[cc]{$\mathsf a_4$}}
\put(70,120){\line(1,-1){20}} \put(76,132){\makebox(0,0)[cc]{$\mathsf a_1$}}
\put(70,120){\line(1,0){40}} \put(90,120) {\makebox(0,0)[bc]{$\max(\mathsf a_1+\mathsf a_3,\mathsf a_2+\mathsf a_4)-\mathsf a_0$}}
\put(110,120){\line(-1,1){20}}\put(105,108){\makebox(0,0)[cc]{$\mathsf a_3$}}
\put(110,120){\line(-1,-1){20}}\put(105,132){\makebox(0,0)[cc]{$\mathsf a_2$}}
\put(10,120){\line(1,1){20}} \put(18,108){\makebox(0,0)[cc]{$\mathsf a_4$}}
\put(10,120){\line(1,-1){20}} \put(18,132){\makebox(0,0)[cc]{$\mathsf a_1$}}
\put(50,120){\line(-1,1){20}}\put(42,108){\makebox(0,0)[cc]{$\mathsf a_3$}}
\put(50,120){\line(-1,-1){20}}\put(42,132){\makebox(0,0)[cc]{$\mathsf a_2$}}
\put(30,140){\line(0,-1){40}}\put(32,120) {\makebox(0,0)[cc]{$\mathsf a_0$}}
\put(60,120){\makebox(0,0)[cc]{$\longleftrightarrow$}}
\end{picture}
\caption{Changing of coordinates under a flip of the triangulation. Only the changing coordinates are shown here, the numbers on the other edges remain unchanged.}
\label{dflipl}
\end{figure}
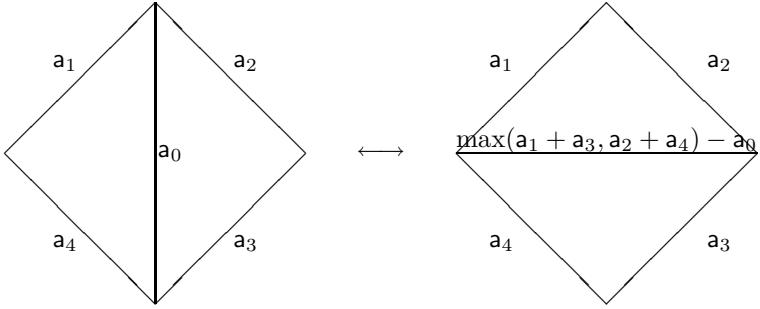

If all edges of the quadrilateral taking part in the flip are different, the coordinates change can be shown by the graphical rule on fig. \ref{dflipl}.

A bounded lamination is integral if and only if all coordinates are integral or half integral and for every triangle with sides $\alpha,\beta,\gamma$ the sum $\mathsf a_\alpha+\mathsf a_\beta+\mathsf a_\gamma$ is integral.

If we multiply a lamination by a positive number, the coordinates get multiplied by the same number. 

A bounded lamination is called {\em even}\index{lamination!even} if all coordinates are integers. This set is denoted by ${\mathsf T}^a_\mathrm{ev}$. One can easily check that  it does not depend on the coordinate system. Observe that an even lamination divided by 2 is not necessarily integral.

The canonical map $\mathsf p\colon \mathsf T^a\rightarrow \mathsf T^x$ is given by $\mathsf x^\alpha=-\sum_\beta\varepsilon^{\alpha\beta}\mathsf a_\beta$.

The action $\mathsf r_\rho$ of the additive group $\mathbb Q$ corresponding to a hole or a cilium $\rho$ amounts to the addition of a number to all coordinates at the edges incident to it.

If $\rho$ is a hole, then the collar map $\mathsf A_\rho$ is given by 
$$\mathsf A_\rho=\max_t(\mathsf a_{\alpha_t}+\mathsf a_{\beta_t}-\mathsf a_{\delta_t}),$$ where $t$ runs over all triangle having $\rho$ as a vertex, $\alpha_t$ and $\beta_t$ are the edges of $t$ incident to $\rho$ and $\delta_t$ is the third edge. 

\paragraph{Real lamination spaces.}

Since the transformation rules for coordinates (\ref{dflipl}) and (\ref{xflipl}) are continuous w.r.t. the standard topology of ${\mathbb Q}^n$, the coordinates define a natural topology on the lamination spaces.  One can now define the spaces of {\em real measured laminations}\index{lamination!real} (resp.  bounded and unbounded) as a completion of the corresponding spaces of rational laminations. These spaces are denoted as $\mathsf T^a(\mathbb R)$ and $\mathsf T^x(\mathbb R)$, respectively. Of course, we have the coordinate systems on these spaces automatically.

  Observe that to define real measured laminations it is not enough just to replace rational numbers by real numbers in the definition of the space of laminations unless the surface $S$ is a disc or a cylinder with one hole. Such definition would not be equivalent to the one above since a sequence of more and more complicated curves with smaller and smaller weights may converge to a real measured lamination which is not presentable by a finite collection of curves.

\section{Teichm\" uller spaces.}

{\em The Teichm\" uller space}\index{Teichm\"uller space} ${\cal T}(S)$ (resp. {\em Moduli space} ${\cal M}(S)$) of a closed surface $S$ is the space of complex structures on $S$ modulo diffeomorphisms isotopic to the identity (resp. modulo all diffeomorphisms). We are going to give two different extensions of this notion to the case of open surfaces with cilia. But before we recall some basic facts about relations between complex structures, constant negative curvature metrics and discrete subgroups of the group $PSL(2,{\mathbb R})$. For more details we recommend the reviews \cite{Abikoff} and \cite{KAG}.

Consider first a vicinity of a boundary component. Topologically it is a cylinder, but as a complex surface it can be isomorphic either to a cylinder or to a punctured disk. The boundary components of the second kind are called {\em punctures}. 

 According to the Poincar\'e uniformisation theorem any hyperbolic complex surface $S$ can be represented as a quotient of the upper half plane $H$ by a discrete subgroup $\Delta$ of its automorphism group (sometimes called the {M\" obius group}) $PSL(2,{\mathbb R})$ of real $2 \times 2$ matrices with unit determinant considered up to the factor $ -1$.  The group $\Delta$ is canonically isomorphic to the fundamental group of the surface $\pi_1(S)$ and it acts freely on $H$. In particular, it implies that the Teichm\"uller space ${\cal T}(S)$ for a surface $S$ without cilia is isomorphic to a connected component in the space of discrete subgroups of the M\"obius group considered up to conjugation. The connected component is singled out by the topology of the quotient of the upper half plane. The group $\Delta$ is called the {\em monodromy group}\index{monodromy group}.
 
Recall that an element $\mathbf g$ of $PSL(2,\mathbb R)$ is called {\em hyperbolic} (resp. {\em parabolic, elliptic}) if it has two distinct real eigenvectors (resp. one real eigenvector, two conjugated complex eigenvectors). Since eigenvectors are stable points of the action of $\mathbf g$ on $\mathbb CP^1$, only parabolic and hyperbolic elements acts on $H$ freely and thus all elements of $\Delta$ must be of these two classes. 

Moreover the quotient of $H$ by the action of the subgroup generated by a parabolic element is a punctured disk. Thus the only elements of $\pi_1(S)$ which can be mapped to parabolic elements are loops surrounding punctures.

On $H$ there exists a unique $PSL(2,{\mathbb R})$-invariant curvature $-1$ Riemannian metric given in the standard coordinates by $\frac{dzd\bar z}{(\Im z)^2}$. It induces a metric on $S$. Since this metric is of negative curvature, any homotopy class of closed curves contains a unique geodesic.

Any element of $\Delta$ has one or two stable points on $\mathbb RP^1$, which is the boundary (absolute) of $H$. The {\em singular set}\index{singular set} of a ciliated surface $S$ is the subset of $\mathbb RP^1$ consisting of all stable points of all elements of the Fuchsian group $\Delta$ and of preimages of all cilia. This set is $\Delta$-invariant by construction and thus its convex hull in $H$ is also $\Delta$-invariant. A {\em convex core}\index{convex core} $S_0$ of the surface $S$ is the quotient of this convex hull by $\Delta$. To describe it geometrically in terms of the original surface $S$ consider the closed geodesics surrounding holes as well as geodesics connecting adjacent cilia on the same boundary components. Cut out of the surface the pieces facing the respective boundary components. These pieces as Riemann manifolds are isomorphic either to the positive quadrant of $H$ or its quotients by a cyclic group. The remaining part of the surface $S$ is just its convex core $S_0$. 

For example, if our surface is a ciliated disk, its convex core is an ideal polygon.

 Homotopy classes of closed curves on a surface are in one-to-one correspondence with the conjugacy classes of its fundamental group. Denote by $\gamma$ an element of $\pi_1(S)$ and by $l(\gamma)$ the length of the corresponding geodesic\index{geodesic length}. Then a simple computation shows that

\begin{equation}
l(\gamma) = \left|\log\frac{\lambda_1}{\lambda_2}\right|,\label{length}
\end{equation}
where $\lambda_1$ and $\lambda_2$ are the eigenvalues of the monodromy along the loop $\gamma$. This number is obviously well defined, i.e., it does not depend on the choices of particular representation of $\pi_1(S)$, of a particular element of $\pi_1(S)$ representing a given loop, and of a particular $2\times 2$ matrix representing an element of $PSL(2,{\mathbb R})$. This formula implies that the length of a geodesic surrounding a puncture is zero.

\subsection{Teichm\" uller space of surfaces with holes ${\cal T}^x(S)$.}\index{Teichm\"uller space!$\mathcal X$-space}

\paragraph{Definition.}
{\em The Teichm\"uller space of a ciliated surface with holes} or {\em the Teichm\"uller $\mathcal X$-space}, denoted by ${\cal T}^x(S)$, is the space of complex structures on $S$  together with the orientation of all holes but the punctures.

It is not {\em a priori} obvious that this space possesses a natural topology in which it is connected, but it will be clarified later. 

For every hole $\rho$ there exists a map $r^\rho:{\cal T}^x(S)\rightarrow \mathbb R_{>0}$ given by the exponent of the length of the geodesics surrounding the hole (resp.  exponent of minus the length of the geodesics surrounding the hole) if the orientation of the hole agrees (resp. disagrees) with the orientation of the surface.

\paragraph{Construction of coordinates.}

Let $\Gamma$ be a triangulation of $S$. For any point of ${\cal T}^x(S)$ we are going to describe a rule for assigning a positive real number to each internal edge of $\Gamma$. The collection of these numbers  will give us a global parameterisation of ${\cal T}^x(S)$.

Every edge of $\Gamma$ connecting two cilia, cilia and a puncture, or two punctures can be made geodesic in a canonical way. In order to make geodesic the edges connecting a cilia or a puncture with a hole first draw a geodesic around each hole and choose an arbitrary point on it. Then deform the edge to a geodesic ray connecting the cilia and the chosen point on the geodesic. Observe that the orientation of the hole induces an orientation of the geodesic. Now move the chosen point along the geodesic in the direction given by this orientation. The limiting position of the ray is just the desired one. By the construction it is spiralling around the hole. The analogous procedure can be done for the edges connecting two holes giving a geodesic spiralling around one closed geodesic at one end and around another one at another end. 

We end up by a triangulation of a part (actually of the convex core $S_0$) of the surface $S$ into ideal triangles. 

 Now lift the triangulation to the upper half plane and consider an edge together with two adjacent triangles forming an ideal quadrilateral. The quadruple of points of this quadrilateral have one invariant under the action of

\noindent
\begin{minipage}{4.7cm}
\begin{picture}(0,120)(0,-20)
\includegraphics[scale=0.25]{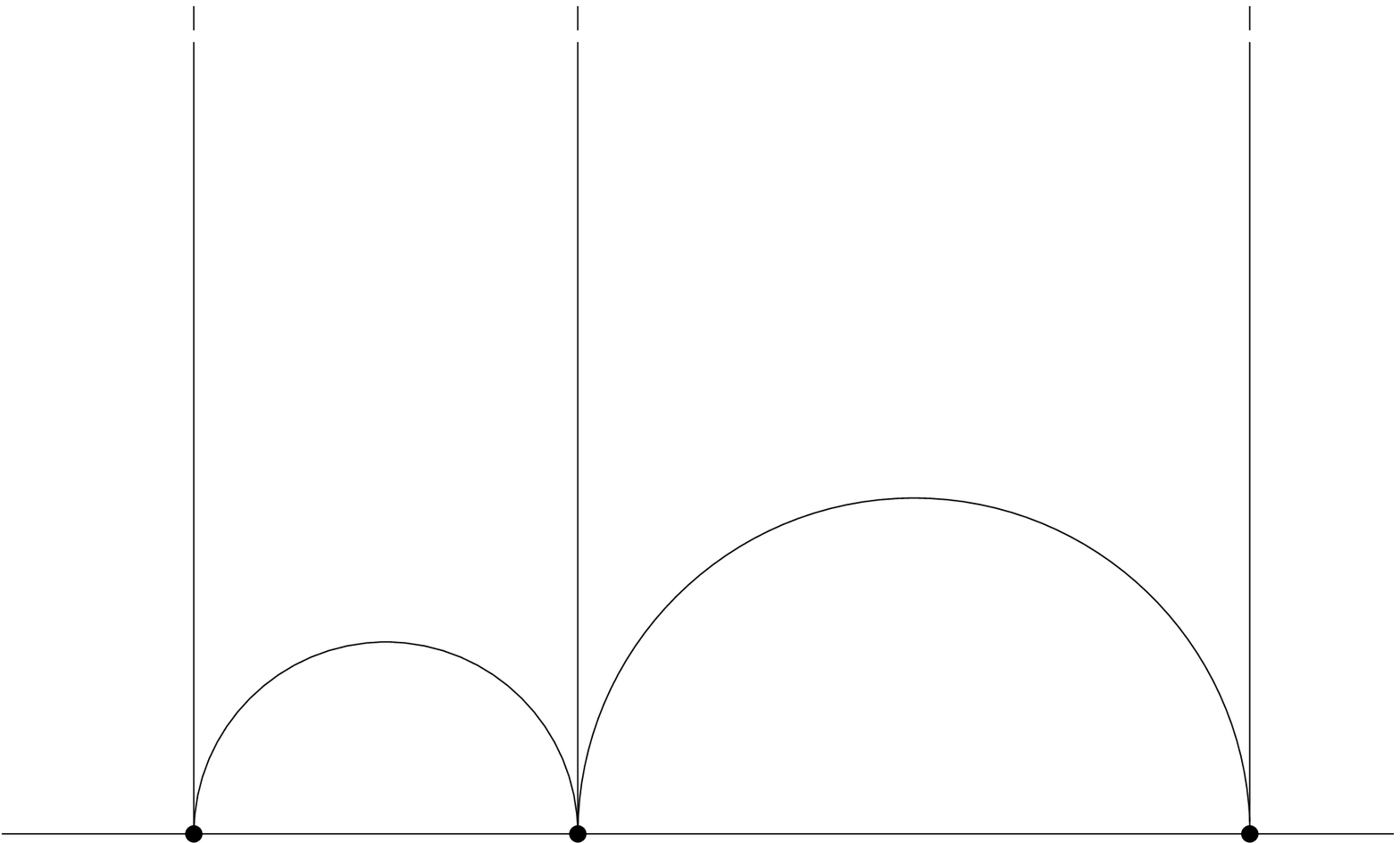}
\put(-62,93){\makebox(0,0)[cc]{$\infty$}}
\put(-14,-6){\makebox(0,0)[cc]{$x$}}
\put(-77,-6){\makebox(0,0)[cc]{$0$}}
\put(-117,-6){\makebox(0,0)[cc]{$-1$}}
\end{picture}
\refstepcounter{figure}\label{cross-ratio}
\noindent Figure \ref{cross-ratio}. Definition of the cross-ratio.
\end{minipage}
\hfill\begin{minipage}{7cm} $PSL(2,\mathbb R)$. For our purpose it is convenient to choose the cross-ratio $x$ by choosing a coordinate on $H$ (sending it to the upper half plane) equal to $0$ and $\infty$ at the ends of the edge and to $-1$ at the third vertex of the quadrilateral (fig. \ref{cross-ratio}). Then the value of the coordinate at the fourth vertex is the desired parameter $x$ associated to the chosen edge. It is an easy exercise to show that there are two ways to choose the coordinate and both give the same value of the cross-ratio. 
\end{minipage}

\paragraph{Reconstruction.}
Our goal now is to construct a ciliated complex surface starting from a triangulation of a ciliated surface with real positive numbers $\{x^\alpha\}$ assigned to internal edges. 

First of all we give a simple receipt how to restore orientations of the boundary components from these data: The orientation of a boundary component corresponding to a hole $\gamma$ is just induced from the orientation of the surface (resp. opposite to the induced one) if the sum $\sum x^\alpha$ over all $\alpha$ incident to $\gamma$ is positive (resp. negative). If the sum is zero, it means that it is not a hole, but a puncture.

Construction of the surface itself can be achieved in two equivalent ways. We shall describe both since one is more transparent from the geometric point of view and the other is useful for practical computations.

{\em Construction by gluing.} We are going to glue the convex core $S_0$ of our surface $S$ out of ideal hyperbolic triangles. The lengths of the sides of ideal triangles are infinite and therefore we can glue two triangles in many ways which differ by shifting one triangle with respect to another along the side. The ways of gluing triangles can be parameterised by the cross-ratios of four vertices of the obtained quadrilateral (considered as points of ${\mathbb R}P^1$) just as defined by fig.~\ref{cross-ratio}.

Now consider the triangulated surface $S$. Replace every triangle of the triangulation by an ideal one and glue them together along the edges just as they are glued in $S$ using numbers assigned to the edges as gluing parameters.

Note, although this is not quite obvious that the resulting surface is not necessarily complete even if the surface has no cilia. In fact it is the complex core $S_0$ with boundary components corresponding to holes removed.
\begin{figure}[h]
\unitlength=1pt
\centerline{\begin{picture}(150,100)(-75,-50)
\thinlines
\put(0,-55){\line( 0, 1){110}}
\put(0,-55){\line( 2, 1){110}}
\put(0,-55){\line(-2, 1){110}}
\put(0, 55){\line( 2,-1){110}}
\put(0, 55){\line(-2,-1){110}}
\thicklines
\put(50,20){\vector(0,-1){ 40}}
\put(50,-20){\vector(-2,1){40}}
\put(10,0){\vector(2,1){40}}
\put(-50,-20){\vector(0,1){ 40}}
\put(-50,20){\vector(2,-1){40}}
\put(-10,0){\vector(-2,-1){40}}
\put(-10,0){\line(1,0){20}}
\put(-50,20){\line(-1,2){10}}
\put(-50,-20){\line(-1,-2){10}}
\put(50,20){\line(1,2){10}}
\put(50,-20){\line(1,-2){10}}
\put(-5,14){\makebox(0,0)[cb]{$\alpha$}}
\put(0,-1){\makebox(0,0)[ct]{$B(x^\alpha$)}}
\put(-30,13){\makebox(0,0)[cb]{$I$}}
\put(-30,-12){\makebox(0,0)[ct]{$I$}}
\put(-50,0){\makebox(0,0)[rc]{$I$}}
\put(30,13){\makebox(0,0)[cb]{$I$}}
\put(30,-12){\makebox(0,0)[ct]{$I$}}
\put(52,0){\makebox(0,0)[lc]{$I$}}
\end{picture}}
\caption{Graph for computing monodromy for $\mathcal T^x(S)$. The original triangulation is shown by the thin lines.}\label{x-monodromy}
\end{figure}
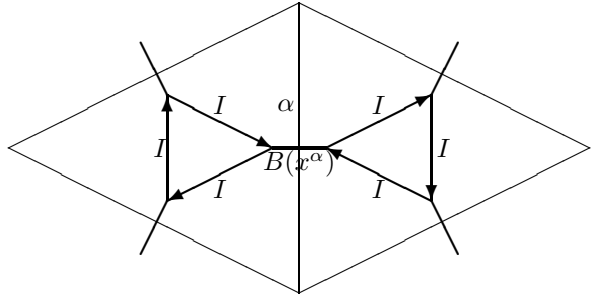 

{\em Construction of the Fuchsian group}. We are now going to construct a discrete monodromy subgroup $\Delta$ of  $PSL(2,{\mathbb R})$ starting from a triangulation of $S$ with positive real numbers on edges. Construct first a graph out of the triangulation in the following way (fig. \ref{x-monodromy}). Draw a small edge of the graph transversal to every edge of the triangulation and for every triangle connect the three ends of edges in it pairwise. Orient the edges of the arising triangle in the clockwise direction. Now assign to each of these edges the matrix $I = \left(\begin{array}{cc} 1 & 1 \\ -1 & 0 \end{array}\right)$. Assign to the remaining edges  the matrix
$B(x^\alpha) = \left(\begin{array}{cc} 0 &
(x^\alpha)^{1/2} \\ -(x^\alpha)^{-1/2} & 0 \end{array} \right)$, where $\alpha$ is the edge of the triangulation it intersects. Now for any oriented path on this graph we can associate a matrix by multiplying consecutively all matrices we meet along it, taking $I^{-1}$ instead of $I$ each time the orientation of the path disagrees with the orientation of the edge. (The orientations of the edges transversal to the edges of the triangulation are not to be taken into account, since $B(x)^2=-1$ and therefore $B(x)$ coincides with its inverse in the group $PSL(2,{\mathbb R})$.) In particular, if we take closed paths starting form a fixed vertex of the graph, we get a homomorphism of the fundamental group of $\Gamma$ to the group $PSL(2,{\mathbb R})$. The image of this homomorphism is the desired group $\Delta$.

 The fact that these two constructions are inverse to the above construction of coordinates is almost obvious, especially for the first one.  The only note we would like to make here is to show where the matrices $I$ and $B(x)$ came from.  Consider two ideal triangles on the upper half plane with vertices at the points ${-1,0,\infty}$ and ${x,\infty,0}$, respectively. Then the M\"obius transform which permutes the vertices of the first triangle is given by the matrix $I$, and the one which maps one triangle to another (respecting the order of vertices given two lines above) is given by $B(x)$.

\begin{figure}[h]
\unitlength 1mm
\begin{picture}(10,50)(4,95)
\put(70,120){\line(1, 1){20}} \put(76,108){\makebox(0,0)[cc]{$x^4(1+(x^0)^{-1})^{-1}$}}
\put(70,120){\line(1,-1){20}} \put(72,132){\makebox(0,0)[cc]{$x^1(1+x^0)$}}
\put(70,120){\line(1,0){40}} \put(90,120) {\makebox(0,0)[bc]{$(x^0)^{-1}$}}
\put(110,120){\line(-1,1){20}}\put(108,108){\makebox(0,0)[cc]{$x^3(1+x^0)$}}
\put(110,120){\line(-1,-1){20}}\put(105,132){\makebox(0,0)[cc]{$x^2(1+(x^0)^{-1})^{-1}$}}
\put(10,120){\line(1,1){20}} \put(18,108){\makebox(0,0)[cc]{$x^4$}}
\put(10,120){\line(1,-1){20}} \put(18,132){\makebox(0,0)[cc]{$x^1$}}
\put(50,120){\line(-1,1){20}}\put(42,108){\makebox(0,0)[cc]{$x^3$}}
\put(50,120){\line(-1,-1){20}}\put(42,132){\makebox(0,0)[cc]{$x^2$}}
\put(30,140){\line(0,-1){40}}\put(32,120) {\makebox(0,0)[cc]{$x^0$}}
\put(60,120){\makebox(0,0)[cc]{$\longleftrightarrow$}}
\end{picture}
\caption{Changing of coordinates under a flip of the triangulation. Only the changing coordinates are shown here, the numbers on the other edges remain unchanged.}
\label{xflip}
\end{figure}
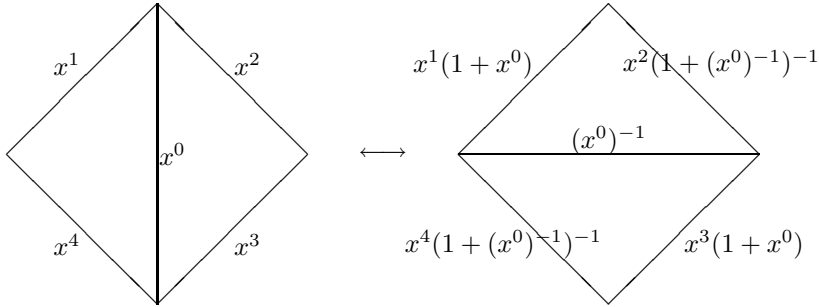

\paragraph{Properties of the coordinates.}
The constructed coordinates on the Teichm\"uller space correspond to a particular choice of the triangulation.  The rule of coordinate change for a flip is given by the formula
\begin{equation}\label{x-flip}
x'^{\beta'}=\left\{\begin{array}{lcl}
(x^\alpha)^{-1}&\mbox{ if }\beta=\alpha\\
x^\beta(1+x^\alpha)^{\varepsilon^{\alpha\beta}} &\mbox{ if }\varepsilon^{\alpha\beta}\geq 0\\
x^\beta(1+(x^\alpha)^{-1})^{\varepsilon^{\alpha\beta}} &\mbox{ if }\varepsilon^{\alpha\beta}\leq 0
\end{array}
\right.
\end{equation}
If all edges of the quadrilateral taking part in the flip are different, the coordinate change can be shown by the graphical rule fig.  \ref{xflip}.

The canonical map $r^\rho\colon \mathcal T^x \rightarrow \mathbb R$ is given by \begin{equation}\label{casimir}
r^\rho= \prod_\alpha x^\alpha
\end{equation}
 where the product is taken over all edges incident to $\rho$ (the edge $\alpha$ is taken twice if both ends of $\alpha$ are incident to $\rho$).

\subsection{Teichm\" uller space of decorated surfaces ${\cal T}^a(S)$.}\index{Teichm\"uller space!$\mathcal A$-space}
 In this section we are going to reproduce some results of Penner \cite{Penner}.

 Before giving a definition of the Teichm\"uller space ${\cal T}^a(S)$, recall what a is a {\em horocycle}\index{horocycle}.  A horocycle in the upper half plane $H$ at $A \in \mathbb RP^1$ is a limit of a circle passing through a fixed point on $H$ when its hyperbolic centre tends to $A$. Equivalently a horocycle at $A$ is an orbit of a subgroup of $PSL(2,\mathbb R)$ stabilising $A$. A horocycle at the point $A \in \mathbb R \subset \mathbb RP^1$ is a circle belonging to the upper half plane tangent to the real axis at $A$. A horocycle at $\infty  \in \mathbb RP^1$ is a horizontal line. Points of one horocycle are equidistant from another horocycle based at the same point. The set of all horocyles is isomorphic to the space $(\mathbb R^2-\{(0,0)\})/\pm 1$. This isomorphism is canonical s equivariant with respect to the $PSL(2)$ action. This isomorphism sends the vector $(x,y)\in \mathbb R^2$ into a horocycle about the point $x/y$ and of euclidean diameter $y^{-2}$.


For a punctured Riemann surface $S = H/\Delta$, where $\Delta$ is a discrete subgroup of $PSL(2,\mathbb R)$, a horocycle at a cilium or a puncture $A$ is a $\Delta$-invariant set of horocycles at each preimage of $A$ on $\mathbb RP^1$.  By definition the area of a horocycle at the puncture is the area inside the horocycle of the quotient of $H$ by the subgroup of $\Delta$ stabilising the basepoint.
 If the horocycle is small enough, its image on $S$ is a small circle surrounding a puncture or tangent to the boundary at the cilium and orthogonal to any geodesics coming out of this puncture or cilium. In this case its area coincides with the actual area inside the image. However the projection of a general horocycle to the surface may have a relatively complicated topology and its area inside the image differs from the area of the horocycle. 

We say that two horocycles on the ciliated surface $S$ are {\em tangent} along the path $\gamma$ connecting their base points if the corresponding horocycles at the ends of the lift of the path to the universal covering $H$ are tangent.  

\paragraph{Definition.}
  A {\em decorated ciliated Riemann surface}\index{surface!decorated} is a ciliated complex surface  such that every hole is a puncture and with a horocycle chosen at each puncture and cilium. The Teichm\" uller space of decorated ciliated surfaces is called the {\em Teichm\" uller space of decorated surfaces} or the {\em Teichm\"uller $\mathcal A$-space}\index{Teichm\"uller space!of decorated surfaces} and is denoted by $\mathcal T^a(S)$.

The space $\mathcal T^a(S)$ has a canonical subspace $\mathcal T^a_0(S)$ consisting of decorated surfaces such that for every path connecting two cilia and retractable to a segment of the boundary without cilia, the corresponding horocycles are tangent.

Forgetting the horocycles we get a canonical map $p\colon\mathcal T^a \rightarrow \mathcal T^x$.

For every cilium and for every hole $\rho$ there is an action of the multiplicative group of positive numbers $r_\rho \colon \mathbb R_{> 0}\times \mathcal T^a\rightarrow \mathcal T^a$ given by changing the size of the corresponding horocycle.

For every hole $\rho$ there is an area map $A_\rho \colon \mathcal T^a \rightarrow \mathbb R_{> 0}$ given by the area inside the corresponding horocycle. 

\paragraph{Construction of coordinates.}
Let $\Gamma$ be a triangulation of $S$. For any point of ${\cal T}^a(S)$, we are going to assign a positive real number to each edge of $\Gamma$. The collection of these numbers  will give us a global parameterisation of $\mathcal T^a(S)$.

\begin{figure}[h]
\includegraphics[scale=0.5]{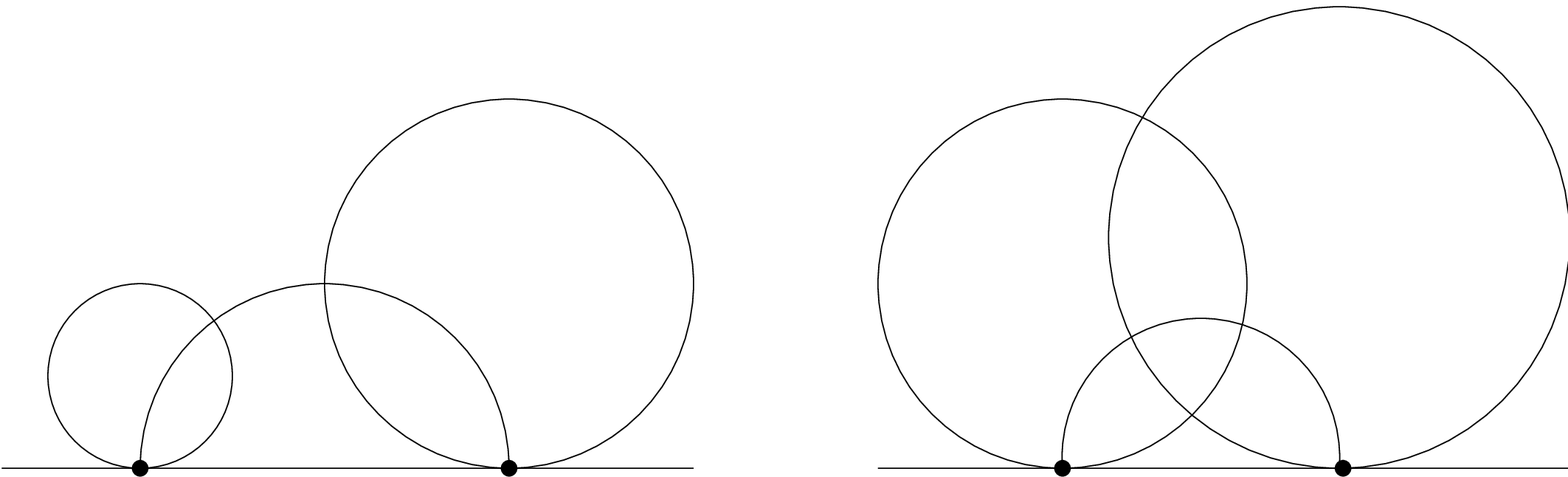}
\put(-74,-6){\makebox(0,0)[cc]{$a_\alpha = e^{-l/2}$}}
\put(-77,35){\makebox(0,0)[cc]{$l$}}
\put(-255,40){\makebox(0,0)[cc]{$l$}}
\put(-245,-6){\makebox(0,0)[cc]{$a_\alpha = e^{l/2}$}}
\caption{Penner coordinates.}
\label{copenner}
\end{figure}

Consider an edge $\alpha$ of the triangulation, deform it into a geodesic and lift it to the universal covering $H$. There are the horocycles attached to the ends of the edges. The number $a_\alpha$ which we assign to the edge is the exponent of a half of the length (resp. of minus the length) of the segment of the lifted edge between the intersection points with the horocycles if the horocycles do not intersect (resp. if the horocycles do intersect), as shown on fig. \ref{copenner}. 

This number can be easily computed algebraically using the correspondence between the horocyles and vectors in $\mathbb R^2$. If $(x_1,y_1)$ and $(x_2,y_2)$ are vectors, corresponding to the horocycles, then $a_\alpha=\left|\det\left(\begin{array}{cc}x_1&x_2\\y_1&y_2\end{array}\right)\right|$.

\paragraph{Reconstruction} is quite similar to that of holed surfaces. There is a canonical mapping $p\colon\mathcal T^a(S) \rightarrow \mathcal T^x(S)$ which just forgets the horocycles and is given explicitly in coordinates by (\ref{tembed}). Therefore to reconstruct the surface itself we can apply the reconstruction procedure for ${\cal T}^x(S)$. To reconstruct the horocycles consider an ideal triangle which we have used to glue the surface. On each edge we have a length of the corresponding geodesics between the horocycles. It allows us to restore unambiguously the points of intersection of the horocycles with the edges.

Observe that given two points $A$ and $B$ on $\mathbb RP^1$ and a horocycle about $A$ one can associate a canonical coordinate $z$ on the upper half plane $H$, such that $z(A)=\infty$, $z(B)=0$ and the horocycle is given by the line $\Im z =1$. Therefore given a lift of an oriented edge of the triangulation to the upper half plane we get a canonical coordinate on the latter taking for $A$ the beginning of this edge and for $B$ its end. Changing the orientation of the edge results in a M\"obius transformation of the coordinate with the matrix $D(a_\alpha)=\left(\begin{array}{cc} 0 &a_\alpha \\ -a_\alpha^{-1} & 0 \end{array} \right)$ and passing from the side $\alpha$ of a triangle with the edges $\alpha, \beta, \gamma$ to the edge $\beta$ with orientations towards the same vertex results in the M\"obius transformation with the matrix $F(\frac{a_\gamma}{a_\alpha a_\beta}) = \left(\begin{array}{cc} 1 & 0 \\  \frac{a_\gamma}{a_\alpha a_\beta}& 1 \end{array} \right)$.

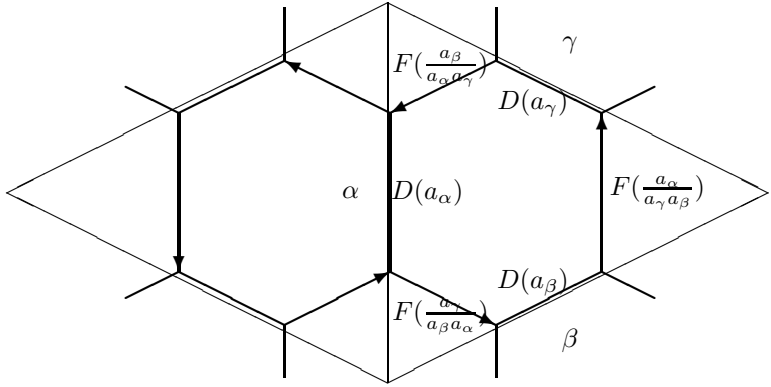
\begin{figure}
\centerline{\begin{picture}(150,150)(-75,-75)
\thicklines
\put(0,-30){\line( 0,1){ 60}}
\put(0,-30){\vector(2, -1){40}}
\put(-40,-50){\vector(2, 1){40}}
\put(40,50){\vector( -2, -1){40}}
\put(0,30){\vector(-2,1){ 40}}
\put(40,-50){\line(2,1){ 40}}
\put(40,50){\line(2,-1){ 40}}
\put(-40,50){\line(-2,-1){ 40}}
\put(-40,-50){\line(-2,1){ 40}}
\put(-80,30){\vector( 0,-1){ 60}}
\put(80,-30){\vector( 0,1){ 60}}
\thinlines
\put(-1,-72){\line( 0,1){ 144}}
\put(-1,-72){\line(2,1){ 144}}
\put(-1, 72){\line(2,-1){144}}
\put(-1,-72){\line(-2,1){ 144}}
\put(-1, 72){\line(-2,-1){144}}
\thicklines
\put(40,50){\line(0,1){ 20}}
\put(40,-50){\line(0,-1){ 20}}
\put(-40,50){\line(0,1){ 20}}
\put(-40,-50){\line(0,-1){ 20}}
\put(-80,-30){\line( -2,-1){ 20}}
\put(80,-30){\line( 2,-1){ 20}}
\put(80,30){\line( 2,1){ 20}}
\put(-80,30){\line( -2,1){ 20}}
\put(-15,0){\makebox(0,0)[cc]{$\alpha$}}
\put(68,56){\makebox(0,0)[cc]{$\gamma$}}
\put(68,-56){\makebox(0,0)[cc]{$\beta$}}
\put(14,0){\makebox(0,0)[cc]{$D(a_\alpha)$}}
\put(54,34){\makebox(0,0)[cc]{$D(a_\gamma)$}}
\put(54,-34){\makebox(0,0)[cc]{$D(a_\beta)$}}
\put(83,0){\makebox(0,0)[lc]{$F(\frac{a_\alpha}{a_\gamma a_\beta})$}}
\put(19,-47){\makebox(0,0)[cc]{$F(\frac{a_\gamma}{a_\beta a_\alpha})$}}
\put(19,47){\makebox(0,0)[cc]{$F(\frac{a_\beta}{a_\alpha a_\gamma})$}}
\end{picture}}
\caption{Graph for reconstructing the monodromy for $\mathcal T^a(S)$.}\label{a-monodromy}
\end{figure}

Similarly to the case of the space $\mathcal T^x(S)$ one can formulate it in terms of a graph with matrices on edges. Construct first a graph out of the triangulation in the following way (Fig. \ref{a-monodromy}). Draw a small edge of the graph along every edge of the triangulation and for every triangle connect the six ends of edges to form a hexagon. Orient the edges surrounding every vertex of the triangulation in the clockwise direction. Assign to each edge of the graph parallel to the edge $\alpha$ the matrix $D(a_\alpha)$ and assign to the edge connecting sides $\alpha$ and $\beta$ of the triangle $\alpha\beta\gamma$ the matrix $F(\frac{a_\gamma}{a_\alpha a_\beta})$. Now for any closed path on the surface one can retract this path to the graph and the monodromy operator can be computed by the product of the corresponding matrices on the edges. 

Moreover this graph with matrices allows to compute the signed distance between any two horocycles. Indeed, consider the path connecting two corresponding vertices of the triangulation. Move each endpoint along an edge of the triangulation to the nearest vertex of the graph and then deform it to pass along the edges of the graph. Associate to it the product of the matrices along this path. Then the absolute value of the upper right element of this matrix is equal to the exponent of signed length between the horocycles. The result does not depend on how we have deformed the path since the monodromy along the edges surrounding the vertex is strictly lower triangular. 

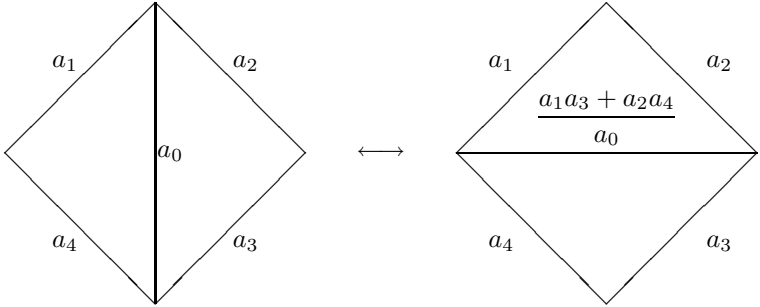
\begin{figure}[h]
\unitlength 1mm
\begin{picture}(100,50)(4,95)
\put(70,120){\line(1, 1){20}} \put(76,108){\makebox(0,0)[cc]{$a_4$}}
\put(70,120){\line(1,-1){20}} \put(76,132){\makebox(0,0)[cc]{$a_1$}}
\put(70,120){\line(1,0){40}} \put(90,121) {\makebox(0,0)[bc]{$\dfrac{a_1a_3+a_2a_4}{a_0}$}}
\put(110,120){\line(-1,1){20}}\put(105,108){\makebox(0,0)[cc]{$a_3$}}
\put(110,120){\line(-1,-1){20}}\put(105,132){\makebox(0,0)[cc]{$a_2$}}
\put(10,120){\line(1,1){20}} \put(18,108){\makebox(0,0)[cc]{$a_4$}}
\put(10,120){\line(1,-1){20}} \put(18,132){\makebox(0,0)[cc]{$a_1$}}
\put(50,120){\line(-1,1){20}}\put(42,108){\makebox(0,0)[cc]{$a_3$}}
\put(50,120){\line(-1,-1){20}}\put(42,132){\makebox(0,0)[cc]{$a_2$}}
\put(30,140){\line(0,-1){40}}\put(32,120) {\makebox(0,0)[cc]{$a_0$}}
\put(60,120){\makebox(0,0)[cc]{$\longleftrightarrow$}}
\end{picture}
\caption{Changing of coordinates under a flip of the triangulation. Only the changing coordinates are shown here, the numbers on the other edges remain unchanged.}
\label{aflip}
\end{figure}

\paragraph{Properties of the coordinates.}
Just like for the three cases considered above, the constructed coordinates correspond to a particular choice of the triangulation.  The rule how the coordinate changes under a flip of an edge $\alpha$ is given by
\begin{equation}\label{a-flip}
a'_{\beta'}=\left\{\begin{array}{lcl}
\dfrac{\prod\limits_{\delta|\varepsilon^{\beta\delta}>0}a_\delta^{\varepsilon^{\beta\delta}}+\prod\limits_{\delta|\varepsilon^{\beta\delta}<0}a_\delta^{-\varepsilon^{\beta\delta}}}{a_\alpha}&\mbox{ if }\beta=\alpha\\
a_\beta &\mbox{ if }\beta\neq\alpha\\
\end{array}
\right.
\end{equation}
If all edges of the quadrilateral taking part in the flip are different the coordinate change can be shown by the graphical rule shown on fig. \ref{aflip}.

The canonical map $p\colon \mathcal T^a\rightarrow \mathcal T^x$ is given by 
\begin{equation}\label{tembed}
x^\alpha=\prod_\beta(a_\beta)^{-\varepsilon^{\alpha\beta}}.
\end{equation}
The action $r_\rho$ of the multiplicative $\mathbb R_{\ge 0}$ corresponding to a hole or a cilium amounts to the multiplication by a number of all coordinates at the edges incident to it.

If $\rho$ is a hole, then the area map $\mathsf A_\rho$ is given by 
$$\mathsf A_\rho=\sum_t(a_{\alpha_t}a_{\beta_t}/a_{\delta_t}),$$
 where $t$ runs over all triangle having $\gamma$ as a vertex, $\alpha_t$ and $\beta_t$ are the edges of $t$ incident to $\gamma$ and $\delta_t$ is the third edge. 

The subspace $\mathcal T^a_0$ is the space where all numbers assigned to external edges are equal to one.

\section{Tropicalisation.}\index{tropicalization}
In this short section we shall state explicitly in which sense the spaces of laminations are ``tropicalisations'' of the corresponding Teichm\"uller spaces. 

Recall that a set $\mathbb F$ with two binary operations: addition $+$ and multiplication $\cdot$ is called a {\em semifield}\index{semifield} if it satisfies the following axioms:
\begin{enumerate}
\item $\mathbb F$ is an Abelian semigroup with respect to addition.
\item $\mathbb F$ is an Abelian group with respect to multiplication. 
\item Addition and multiplication are related by the distributivity law:\\ $a\cdot(b+c)=a\cdot b+a\cdot c$.
\end{enumerate}
The simplest examples of a semifield are just positive real or rational numbers with respect to the usual addition and multiplication. Another class of examples are given by the sets of real, rational and integral numbers with the usual addition used for multiplication and usual maximum used for addition. The corresponding semifields are called {\em real, rational and integral tropical semifield} and are denoted by $\mathbb R^t,\mathbb Q^t$ and $\mathbb Z^t$ repectively. Every algebraic expression which does not contain substraction makes sense in a semifield. Consider for example the polynomial $a^3+3ab^2+2$. In tropical semifields this expression can be rewritten using usual operation as $\max(3\mathsf a,2\mathsf a + \mathsf b,0)$. We say that the second expression is the {\em tropical analogue}\index{tropical analogue} of the first one. 

Observe that tropical semifields can be considered as certain limits of nontropical ones. Namely fix a positive real constant $\epsilon$ and consider the set of real numbers provided by the addition:$(a,b)\mapsto \epsilon\log(e^{a/\epsilon}+e^{b/\epsilon})$ and the multiplication given by the ordinary addition. This semifield is obviously isomorphic to the semifield of positive numbers by the isomorphism $a\mapsto e^{a/\epsilon}$. In the limit $\epsilon \rightarrow 0$ the isomorphism is not longer defined, but the limit of the operations does exist and coincides with the operations in the tropical semifield. 

In particular the tropical analogue of a rational function of $n$ variables $f(a_1,\ldots,a_n)$ which does not contain substraction is given by $\mathsf f(\mathsf a_1,\ldots,\mathsf a_n)=\lim_{\epsilon \rightarrow 0} \epsilon \log f(e^{\mathsf a_1/\epsilon},\ldots,e^{\mathsf a_n/\epsilon})$.

One can easily see by comparison of the sections 3 and 4 that the formulae for the spaces of laminations are tropical analogues of the respective formulae for the Teichm\"uller spaces. For example space of laminations $\mathsf T^x$ has he same set of coordinate systems as the Teichm\"uller space $\mathcal T^x$, while the transition functions between the coordinate system (\ref{xl-flip}) for the former space are given by tropical analogues of the respective functions (\ref{x-flip}) for the latter one. Similarly the lamination spaces $\mathsf T^a$ and $\mathsf T^a_0$ are tropical analogues of the spaces $\mathcal T^a$ and $\mathcal T^a_0$ respectively.

In particular one can see that Teichm\"uller spaces can be canonically compactified by adding the quotients of the respective lamination spaces by action of the multiplicative group of positive real numbers. This is just a verision of the {\em Thurston compactification}\index{Thurston compactification} \cite{Thurston} for ciliated surfaces.

\section{Poisson and degenerate symplectic structures.}\index{Weil-Petersson form}
We are going to give a definition (by explicit formulae) of a Poisson bracket $P_{\mbox{\small \it WP}}$ on the spaces $\mathcal T^x(S)$ and of a degenerate symplectic structure $\omega_{\mbox{\small \it WP}}$ on $\mathcal T^a(S)$. For surfaces $S$ without cilia such structures are well known under the name of Weil-Petersson ones and can be defined as follows.

The spaces of homomorphisms of the fundamental group of a surface $S$ into any reductive group $G$ considered up to conjugation is known to have a canonical Poisson structure (see \cite{FR} and references therein). Its symplectic leaves are numerated by conjugacy classes of elements corresponding to loops surrounding holes. The Teichm\"uller space $\mathcal T^x(S)$ can be mapped to this space (for $G=PSL(2,\mathbb R)$) and if the surface $S$ has no cilia the corresponding map is a local diffeomorphism in a vicinity of a generic point. This map induces the Poisson structure $P_{\mbox{\small \it WP}}$. 

Another way of defining the same Poisson structure was given by William Goldman in \cite{Goldman1}. For any closed loop one can associate the absolute value of the trace of the element of $PSL(2,\mathbb R)$ corresponding to it. Goldman gave an explicit expression for the Poisson bracket between function of this type as a linear combination of functions of the same type.  

The image of $\mathcal T^a(S)$ in $\mathcal T^x(S)$ under the map $p$ corresponds to the representations of $\pi_1(S)$ in $PSL(2,\mathbb R)$ sending all loops surrounding holes  to parabolic elements. It implies that the Poisson structure restricted to the image is nondegenerate. The inverse image under the map $p$ of the corresponding symplectic structure defines a degenerate symplectic structure $\omega_{\mbox{\small \it WP}}$ on $\mathcal T^a(S)$.

The Poisson structure on $\mathcal T^x(S)$ for a ciliated surface $S$ is given by the formula:
\begin{equation}\label{poisson}
P_{\mbox{\small \it WP}}=\sum_{\alpha\beta\in E(\Gamma)}\varepsilon^{\alpha\beta}x^\alpha \frac{\partial}{\partial x^\alpha}\wedge x^\beta\frac{\partial}{\partial x^\beta}.
\end{equation}

The degenerate symplectic structure $\omega_{WP}$ on $\mathcal T^a(S)$ for ciliated surface $S$ is given by the formula:
\begin{equation}\label{symplectic}
\omega_{\mbox{\small \it WP}}=\sum_{\alpha\beta\in E(\Gamma)}\varepsilon^{\alpha\beta}\frac{da_\alpha}{a_\alpha}\wedge \frac{da_\beta}{a_\beta}
\end{equation}

The main property of these Poisson and symplectic structures is that they are independent of the particular triangulation. This property can be checked explicitly by substitution of the transformation laws (\ref{x-flip}) (resp. (\ref{a-flip})) into the expression for $P_{\mbox{\small \it WP}}$ (resp. $\omega_{\mbox{\small \it WP}}$), taking into account the transformation law (\ref{epsilon-flip}) for the matrix $\varepsilon^{\alpha\beta}$.

Another useful property is the compatibility of the bracket with unramified covering of surfaces. If $\tilde S\rightarrow S$ is an unramified covering of degree $N$ then there are obvious maps $i^x:\mathcal T^x(S) \rightarrow \mathcal T^x(\tilde S)$ and $i^a:\mathcal T^a(S) \rightarrow \mathcal T^a(\tilde S)$. Then one can easily see that $(i^a)^*\tilde\omega_{\mbox{\small \it WP}} =N\omega_{\mbox{\small \it WP}}$, where $\tilde\omega_{\mbox{\small \it WP}}$ is the canonical degenerate symplectic form on $\mathcal T^x(\tilde S)$.   The map $i^x$ is not a Poisson map, however it maps symplectic leaves to symplectic leaves. If $i^x:L\rightarrow \tilde L$, where $L$ and $\tilde L$ are the symplectic leaves of $\mathcal T^x(S)$ and  $\mathcal T^x(\tilde S)$ respectively, then $(i^x)^*(\tilde P_{\mbox{\small \it WP}})^{-1}|_{\tilde L}=N(P_{\mbox{\small \it WP}})^{-1}|_ L$.

One can check that if the surface $S$ has no cilia the formula gives the known Poisson and symplectic structures. The simplest way to do it is to compute explicitly the bracket between traces of matrices corresponding to closed loops using both our definition and Goldman's one. Due to the compatibility with coverings it is sufficient to check the coincidence of the brackets for two loops having one intersection point when Goldman's expression is especially simple. We do not give the computation here since it is straightforward and technical.  

For the surface $S$ without cilia the kernel of the symplectic form is spanned by the vector fields corresponding to the action $r_\rho$, where $\rho$ runs through the holes. The kernel of the Poisson structure is spanned by the differentials of the functions $r^\rho$. 

The lamination spaces $\mathsf T^x$ and $\mathsf T^a$ also have the canonical Poisson structure $\mathsf P_{\mbox{\small \it WP}}$ and the degenerate symplectic structure $\mathsf w_{\mbox{\small \it WP}}$, respectively. They are given by 
$$
\mathsf P_{\mbox{\small \it WP}}=\sum_{\alpha\beta\in E(\Gamma)}\varepsilon^{\alpha\beta}\frac{\partial}{\partial \mathsf x^\alpha}\wedge \frac{\partial}{\partial \mathsf x^\beta}.
$$
$$
\mathsf w_{\mbox{\small \it WP}}=\sum_{\alpha\beta\in E(\Gamma)}\varepsilon^{\alpha\beta}d\mathsf a_\alpha\wedge d\mathsf a_\beta
$$

All the properties of these structures are similar to the respective properties of the Poisson and symplectic structures on their Teichm\"uller counterparts.

\section{Canonical pairing.}\label{Pairing}

Certain curves or collections of curves on 2D surface define functions on Teichm\"uller and lamination spaces. Namely:
\begin{enumerate}
\item\label{closed-on-teich} A closed curve $\gamma$ on $S$ defines a function $\ell(\gamma)$ on $\mathcal T^a(S)$ as well as on $\mathcal T^x(S)$ by the length of the geodesic w.r.t. the curvature $-1$ metric in the homotopy class defined by the curve.
\item\label{open-on-a} If a curve $\gamma$ on $S$ connects two points on the boundary it defines a function $\ell(\gamma)$ on $\mathcal T^a(S)$ in the following way. Move the ends of the curve which do not belong to holes, in the counterclockwise direction until they hit the cilia. Then take a geodesic in the same homotopy class and take the exponent of the signed half of the length of the segment between the corresponding horocycles like on fig. \ref{copenner}.
\item\label{open-on-x} If a collection of curves $\mathsf m = \{\gamma_i\}$ with rational weights $\{w_i\}$ has the property that for every segment of the boundary between cilia and for every hole the total weight of the curves hitting it vanishes, then it defines a function on $\mathcal T^x(S)$ in the following way. Move first the ends of curves in the counterclockwise direction to hit the cilia. Lift then the collection of curves to $H$ and choose horocycles at their ends. Then take the weighted sum of signed lengths of geodesic segments between the horocycles $\ell(\mathsf m)=\sum w_i \ell(\gamma_i)$. It is easy to see that this sum does not depend on the choice of the horocycles due to the vanishing property. 
\end{enumerate}

These functions can be considered as functions of two arguments: the first is a curve or a weighted collection of curves and the second is a point of the appropriate Teichm\"uller or lamination space. Our aim now is to extend the first argument to be not just a single curve, but a point of an appropriate space of laminations. The resulting functions of two arguments are called {\em pairings}. 

There are two version or pairings between Teichm\"uller and lamination spaces --- additive and multiplicative. The main feature of the additive one is that it can be extended to real laminations as a continuous function. The multiplicative one is defined for integral laminations only. For every lamination this pairing defines a function on the Teichm\"uller space being a Laurent polynomial in the coordinates with positive integral coefficients.

\subsection{Additive canonical pairing.}\index{pairing!additive}
The additive canonical pairing  is a function of two arguments. One argument is a rational (or, later on, real) lamination, which can be either from $\mathsf T^a_0$ or $\mathsf T^x$, and the other is a point of the opposite type  Teichm\"uller space, $\mathcal T^x$ or $\mathcal T^a_0$, respectively.
$$
\mathcal I\colon \mathsf T^x \times \mathcal T^a_0 \rightarrow \mathbb R, 
$$
$$
\mathcal I\colon \mathcal T^x \times \mathsf T^a_0\rightarrow \mathbb R, 
$$
Abusing notation,  we shall denote all of them by a single letter $\mathcal I$. 

We will also define the {\it intersection pairing}:\index{pairing!intersection}
$$
\mathsf I\colon \mathsf T^x \times \mathsf T^a_0\rightarrow \mathbb R
$$

Sometimes it is useful to consider the functions $\mathcal I$ and $\mathsf I$ as a function of one variable taking another one as a parameter. In this case we shall put this parameter as the lower index. For example $\mathcal I_{\mathsf m^x}(m^a):=\mathbb I(\mathsf m^x,m^a)$.

The additive and the intersection pairings are defined as follows. 
\paragraph{Definition}
\begin{enumerate}
\item Let $\gamma^x \in \mathsf T^x $  be a single closed curve and $m^a\in \mathcal T^a_0$. Then $\mathcal I(\gamma^x,m^a)$ is its length $\ell(\gamma^x)$  (rule \ref{closed-on-teich}).

\item Let $m^x\in \mathcal T^x$ and let $\gamma^a \in \mathsf T^a$ be a single closed curve.  Then $\mathcal I(m^x,\gamma^a)$ is equal to $\pm$ a half of its length $\ell(\gamma^a)$ (rule \ref{closed-on-teich}).  The sign is $+$ unless $\gamma^a$ surrounds a negatively oriented hole.

\item Let $\gamma^x\in \mathsf T^x$ be a curve connecting two points of the boundary of $S$ and $m^a \in \mathcal T^a$. Then $\mathcal I(\gamma^x,m^a)$ is the signed length $\ell(\gamma^x)$ given by the rule \ref{open-on-a}. 

\item Let $\gamma^x\in \mathsf T^x$ and $\gamma^a\in \mathsf T^a$ be two curves. Move the ends of $\gamma^a$ (if any and but those belonging to holes) counterclockwise to hit cilia.  Then $\mathsf I(\gamma^x,\gamma^a)$ is a half of the minimal number of intersection points between $\gamma^x$ and $\gamma^a$.

\item  Let $m^x\in \mathcal T^x$ and let $\mathsf m^a\in \mathsf T^a_0$ be a collection of curves. Then the function $\ell(l^a)$ defined by this rule \ref{open-on-x} is equal to $\mathcal I(m^x,\mathsf m^a)$. (The vanishing property follows from the definition of  $\mathsf T^a_0$.

\item Let $u,v$ be positive rational numbers. Let $\mathsf m^x_1,\mathsf m^x_2\in \mathsf T^x$ be two collections of nonintersecting curves such that no curve from $\mathsf m^x_1$ intersects $\mathsf m^x_2$,  and $m^a\in {\cal T}^a$. Then  $\mathcal I(u \mathsf m^x_1+v \mathsf m^x_2,m^a) = u \mathcal I(\mathsf m^x_1,m^a)+v \mathcal I(\mathsf m^x_2,m^a)$. If $m^x \in \mathcal T^x$ and  $\mathsf m^a_1,\mathsf m^a_2\in \mathsf T^a_0$ be nonintersecting collections of curves such that no curve from $\mathsf m^a_1$ intersects $\mathsf m^a_2$ then $\mathcal I(m^x, u \mathsf m^a_1+v \mathsf m^a_2) = u \mathcal I(m^x,\mathsf m^a_1)+v \mathcal I(m^x,\mathsf m^a_2)$.
\end{enumerate} 

\paragraph{Properties of the additive canonical pairing.}
\begin{enumerate}

\item\label{ax-continuity} The additive canonical pairings  $\mathcal I$ as well as the intersection pairing $\mathsf I$ are continuous functions.

\item\label{symmetry-add} 
$\mathcal I(p(m^a),\mathsf m^a)=\mathcal I(\mathsf p(\mathsf m^a),m^a)$, where $\mathsf m^a \in {\mathsf T}^a_0$ and $m^a\in {\mathcal T}^a_0$.

\item\label{explicit} If a point $\mathsf m^x$ of the space ${\mathsf T}^x$ has positive coordinates $\mathsf x^1,\ldots,\mathsf x^n$ and a point $m^a$ of the Teichm\"uller space ${\mathcal T}^a_0$ has coordinates $a_1,\ldots,a_n$, then $\mathcal I(\mathsf m^x, m^a)= \sum_\alpha \mathsf x^\alpha \log a_\alpha$.

\item\label{limit-add} Let $m^x \in {\mathcal T}^x$ be a point of a Teichm\"uller space ${\mathcal T}^x$ with coordinates $x_1, \ldots, x_n$. Let $C$ be a positive real number. Denote by $(m^x)^C\in {\mathcal T}^x$ a point with coordinates $(x_1)^C,\ldots,(x_n)^C$ and by $\mathsf m^x \in {\mathsf T}^x$ a lamination with coordinates $\mathsf x^1=\log x^1,\ldots,\mathsf x^n=\log x^n$. Then for any lamination $\mathsf m^a\in {\mathsf T}^a$ we have
$$\lim\limits_{C \rightarrow \infty}\mathcal I((m^x)^C,\mathsf m^a)/C=\mathsf I(\mathsf m^x,\mathsf m^a)$$.

\item Analogously let $m^a\in {\mathcal T}^a_0$ be a point of the Teichm\"uller space with coordinates $a_1,\ldots,a_n$. Let $C$ be a positive real number. Denote by $(m^a)^C \in {\mathcal T}^a_0$ the point in the Teichm\"uller space with coordinates $(a_1)^C,\ldots,(a_n)^C$ and by $\mathsf m^a\in {\mathsf T}^a_0$ the lamination with coordinates $\mathsf a_1=\log a_1,\ldots,\mathsf a_n=\log a_n$. Then for any lamination $\mathsf m^x\in {\mathsf T}^x$ we have
$$\lim\limits_{C \rightarrow \infty}I(\mathsf m^x,(m^a)^C)/C=\mathsf I(\mathsf m^x, \mathsf m^a)$$
\end{enumerate}

For closed surfaces there exists a canonical pairing $\mathcal I\colon \mathcal T^x \times \mathcal T^a_0 \rightarrow \mathbb R$ defined by Francis Bonahon \cite{Bonahon}, such that the pairings defined here are its limits. Unfortunately we don't know any explicit construction of this pairing. 

\subsection{Multiplicative canonical pairing.}\index{pairing!multiplicative} 
The multiplicative canonical pairing  is a function of two arguments. One argument is an integral lamination which can be either from $\mathsf T^a_0(\mathbb Z)$ or $\mathsf T^x(\mathbb Z)$, and the other is a point of the opposite type  Teichm\"uller space, $\mathcal T^x$ or $\mathcal T^a_0$, respectively.
$$
\mathbb I\colon \mathsf T^x(\mathbb Z) \times \mathcal T^a_0 \rightarrow \mathbb R_{>0}, 
$$
$$
\mathbb I\colon \mathcal T^x \times \mathsf T^a_0(\mathbb Z)\rightarrow \mathbb R_{>0}, 
$$

Abusing notations, we shall denote all of them by a single letter $\mathbb I$. 
Sometimes it is useful to consider the function $\mathbb I$ as a function of one variable taking another one as a parameter. In this case we shall put this parameter as the lower index. For example $\mathbb I_{\mathsf m^x}(m^a):=\mathbb I(\mathsf m^x,m^a)$.

The multiplicative pairing is defined by the following properties: 
\paragraph{Definition}
\begin{enumerate}
\item Let  $\gamma^x \in \mathsf T^x $  be a single closed curve with weight $k$ and $m^a\in \mathcal T^a$. Then $\mathbb I(\gamma^x,m^a)$ is the absolute value of the trace of the monodromy of $(\gamma^x)^k$.

\item Let $m^x\in \mathbb T^x$ and let $\gamma^a \in \mathsf T^a$ be a single closed curve which is not retractable to a hole. Then $\mathbb I(\gamma^x,m^a)$ is the absolute value of the trace of the monodromy of $(\gamma^x)^k$.

\item Let $m^x\in \mathcal T^x$ and let $\gamma^a \in \mathsf T^a$ be a single closed curve which is retractable to a hole and the hole is oriented positively (resp. negatively). Then $\mathcal I(\gamma^x,m^a)$ is equal to the absolute value of the largest (resp. smallest) eigenvalue of the monodromy of $(\gamma^x)^k$. 

\item Let $\gamma^x\in \mathsf T^x$ be a curve with weight $k$ connecting two points of the boundary of $S$ and $m^a \in \mathcal T^a$. Then $\mathbb I(\gamma^x,m^a)=\exp \mathcal I(\gamma^x,m^a)$. 

\item  Let $m^x\in \mathcal T^x$ and let $\mathsf m^a\in \mathsf T^a_0$ be a collection of open curves with integral weights. Then $\mathbb I(m^x,\mathsf m^a)= \exp\mathcal I(m^x,\mathsf m^a)$. 

\item Let $u,v$ be positive integral numbers. Let $\mathsf m^x_1,\mathsf m^x_2\in \mathsf T^x$ be two collections of nonintersecting curves such that no curve from $\mathsf m^x_1$ intersects or coincides with a curve from $\mathsf m^x_2$, and $m^a\in {\cal T}^a$. Then  $\mathbb I(u \mathsf m^x_1+v \mathsf m^x_2,m^a) = (\mathbb I(\mathsf m^x_1,m^a))^u(\mathbb I(\mathsf m^x_2,m^a))^v$. 

If $m^x \in \mathcal T^x$ and  $\mathsf m^a_1,\mathsf m^a_2\in \mathsf T^a_0$ be nonintersecting collections of curves such that no curve from $\mathsf m^x_1$ intersects or coincides with a curve from $\mathsf m^x_2$. Then $\mathbb I(m^x, u \mathsf m^a_1+v \mathsf m^a_2) = (\mathbb I(m^x,\mathsf m^a_1))^u(\mathbb I(m^x,\mathsf m^a_2))^v$.
\end{enumerate} 

\paragraph{Properties of the multiplicative canonical pairing.}
\begin{enumerate}
\item\label{x-positivity} Let $\mathsf a_1,\ldots,\mathsf a_n$ be coordinates of a lamination $\mathsf m^a \in \mathsf T^a_0(S,\mathbb Z)$. Then the function $\mathbb I_{\mathsf m^a}$ on $\mathcal T^x$ is a Laurent polynomial with positive integral coefficients of $(x^1)^{1/2},\ldots,(x^n)^{1/2}$, where $x^1,\ldots,x^n$ are the coordinates on $\mathcal T^x$. Its highest term is equal to $\prod_\alpha (x^\alpha)^{\mathsf a_\alpha}$ and the lowest one is $\prod_\alpha (x^\alpha)^{-\mathsf a_\alpha}$. Moreover $\mathbb I_{\mathsf m^a}$ is a Laurent polynomial of $x^1,\ldots,x^n$ if and only if the lamination $\mathsf m^a$ is even.  

\item\label{a-positivity}\index{Laurent property} Let $\mathsf x^1,\ldots,\mathsf x^n$ be coordinates of a lamination $\mathsf m^x \in \mathsf T^x(S,\mathbb Z)$. Then the function $\mathbb I_{\mathsf m^a}$ on $\mathcal T^x$ is a Laurent polynomial with positive integral coefficients of the coordinates $a^1,\ldots,a^n$  on $\mathcal T^a_0$.  

\item\label{rel-to-additive} Let $\mathsf m^x \in \mathsf T^x(S,\mathbb Z)$ be an integral $\mathcal X$-lamination and let $m^a\in \mathcal T^a_0$ be a point of Teichm\"uller $\mathcal A$-space. Then $\lim_{C \rightarrow \infty}\log (\mathbb I(C\mathsf m^x,m^a))/C=\mathcal I(\mathsf m^x,m^a)$. Analogously let $m^x \in \mathcal T^x$ be a point of the Teichm\"uller $\mathcal X$-space and let $\mathsf m^a\in \mathsf T^a_0$ be an $\mathcal A$-lamination. Then \\$\lim_{C \rightarrow \infty}\log (\mathbb I(Cm^x,\mathsf m^a))/C=\mathcal I(\mathsf m^x,m^a)$.

\item\label{limit-mult} Let $\mathsf m^a\in \mathsf T^a_0$ be an integral bounded lamination. Then the function $\mathsf I_{\mathsf m^a}$ on $\mathsf T^x$ is the tropical analogue of the function $\mathbb I_{\mathsf m^a}$ on $\mathcal T^x$ in any coordinate system. Analogously let $\mathsf m^x\in \mathsf T^x$ be an integral unbounded lamination. Then the function $\mathsf I_{\mathsf m^x}$ on $\mathsf T^a_0$ is the tropical analogue of the function $\mathbb I_{\mathsf m^x}$ on $\mathcal T^a$ in any coordinate system.

\item\label{particular-case} The coordinate functions on the space $\mathcal T^a$ as well as  $\mathcal T^a_0$ are particular cases of the pairing, namely the value of the coordinate function $a_\alpha$ at the point $m^a$ of the Teichm\"uller space $\mathcal T^a$ is $\mathbb(\mathsf m^x,m^a)$, where $\mathsf m^x$ is the edge $\alpha$ of the triangulation with weight $1$. 
\end{enumerate}

Conjecturally the converses to the properties 1 and 2 are also true. Namely, theonly functions on $\mathcal T^x$ (resp. $\mathcal T^a_0$) given by a Laurent polynomial with positive integral coefficients in every coordinate chart and undecomposable into a sum of two such polynomials are given by the multiplicative canonical pairing.

\subsection{Proofs of main properties.}

Here we are going to give main ideas of the proofs of the most important properties of pairings: the continuity of the additive pairing (property \ref{ax-continuity}) and the properties \ref{x-positivity} and \ref{a-positivity} of the multiplicative pairing. Other properties are either obvious or follow from the proof of these three.

\paragraph{Continuity}

Let us first prove the continuity of the pairing between ${\mathsf T}^a$ and the ${\mathcal T}^x$. To do this it suffices to prove it for laminations without curves with negative weights. Indeed, if we add such a curve to a lamination the length obviously changes continuously. We are going to show that the length of integral lamination is a convex function of its coordinates, i.e. that

\begin{equation} \label{ineq}
I(m^x,\mathsf m^a_1) + I(m^x,\mathsf m^a_2) \leq I(m^x,\mathsf m^a_{1+2}) 
\end{equation}
were by $\mathsf m^a_{1+2}$ we have denoted a lamination with coordinates being sums of the respective coordinates of $\mathsf m^a_1$ and $\mathsf m^a_2$. Taking into account the homogeneity of $I$, one sees that the inequality (\ref{ineq}) holds for all rational laminations and therefore can be extended by continuity to all real laminations.

Let us prove the inequality (\ref{ineq}).  Draw both laminations $\mathsf m^a_1$ and $\mathsf m^a_2$ on the surface and deform them to be geodesic. These laminations in general intersect each other in finite number of points. Construct now a new collection of curves  

\vspace{2pt}
\noindent
\begin{minipage}{4.2cm}
\begin{picture}(0,120)(-20,-20)
\includegraphics[scale=0.25]{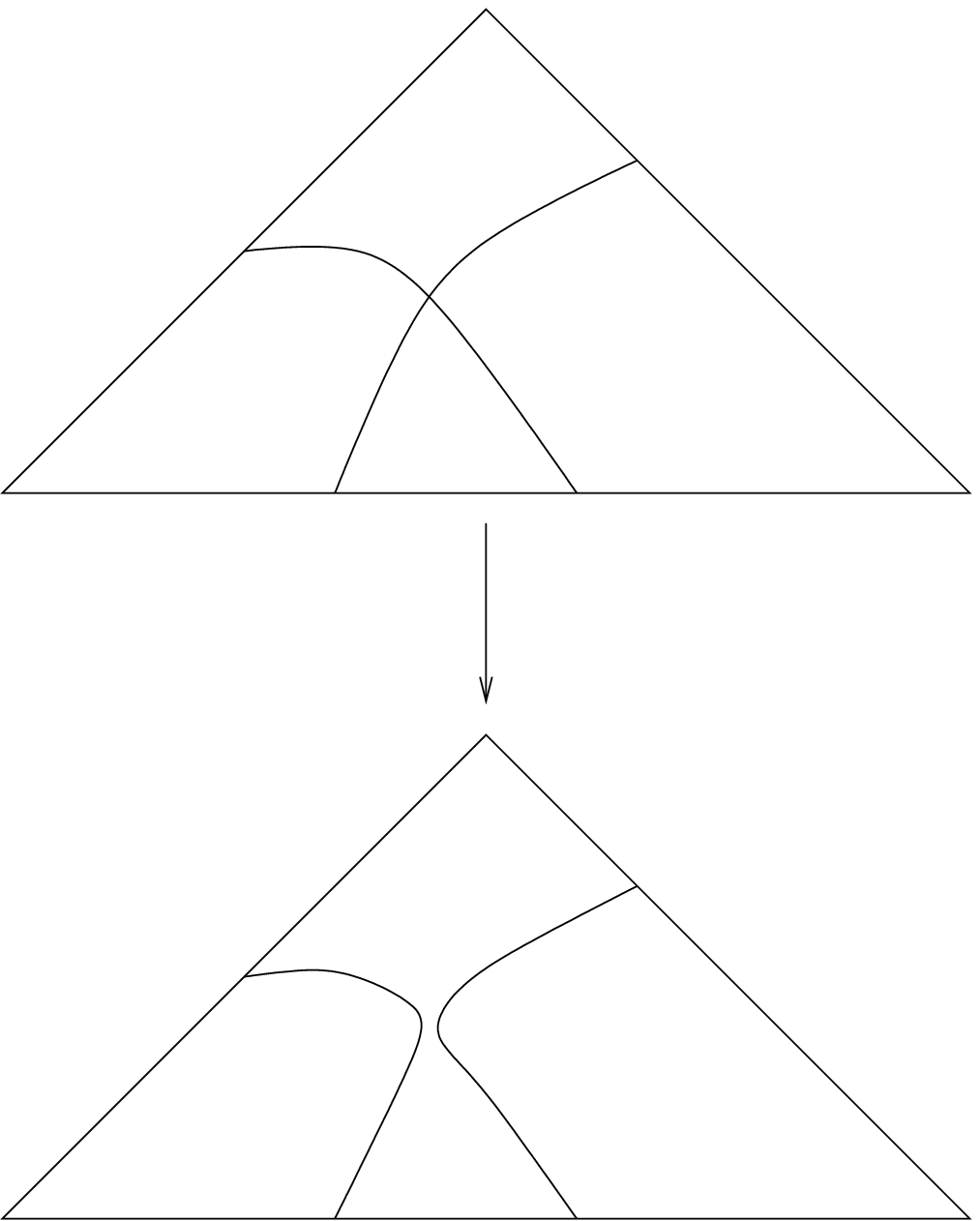}
\put(-43,93){\makebox(0,0)[cc]{$A$}}
\put(-78,0){\makebox(0,0)[cc]{$B$}}
\put(6,0){\makebox(0,0)[cc]{$C$}}
\put(-43,40){\makebox(0,0)[cc]{$A$}}
\put(-78,53){\makebox(0,0)[cc]{$B$}}
\put(6,53){\makebox(0,0)[cc]{$C$}}
\end{picture}
\refstepcounter{figure}\label{continuity}
\noindent Figure \ref{continuity}. Rearranging the intersection points.
\end{minipage}
\hfill\begin{minipage}{7cm} out of the union of $\mathsf m^a_1$ and $\mathsf m^a_2$ in
the following way. For every intersection points consider the triangle $ABC$ of the triangulation it belongs to. Without loss of generality we may assume, that the segments of the curves containing the intersection points connect the sides $AB$ with $BC$ and $AC$ with $BC$, respectively. Then remove the neighbourhood of the intersection point and replace it by a union of two nonintersecting arcs in such a way that the segments of the resulting curves still connect $AB$ with $BC$ and $AC$ with $BC$ (Fig. \ref{continuity}).  The resulting collection of curves  considered
\end{minipage}

\vspace{2pt}
\noindent   as a lamination is just $\mathsf m^a_{1+2}$. The total length of the constructed curves can be made as close to $I(m^x,\mathsf m^a_1) + I(m^x,\mathsf m^a_2)$ as one desires. However it is not geodesic any longer. When we deform this collection to the geodesic one its length only decreases, thus giving the desired inequality~\ref{ineq}.

\paragraph{Laurent property and positivity.}
Now we are going to show that the functions $\mathbb I_{\mathsf m^a}$ and $\mathbb I_{\mathsf m^x}$ are given by Laurent polynomials with positive integral coefficients in any coordinate charts. We shall show it for laminations representable by a single curve and then extend the result to an arbitrary lamination by multiplicativity. 

Consider first a single closed curve with weight $k$ as an element of the space $\mathsf T^a$. If $M$ is the monodromy operator around this curve, then by definition $\mathbb I_{\mathsf m^a}=tr (M^k)$. According to the construction of the Fuchsian group from section 3.1 the operator $M$ is given by the expression $$M= B(x^{i_1})I^{\pm 1}B(x^{i_2})\ldots I^{\pm 1}B(x^{i_p})I^{\pm 1},$$ where $i_1,\ldots,i_p$ is the sequence of edges of the triangulation the path intersects and the signs for $I^{\pm 1}$ are chosen according to whether the path turns left or right passing through a triangle. Multiplying every $I^{\pm 1}$  by $J=\left(\begin{array}{cc}0&1\\-1&0\end{array}\right)$ from the left and every $B$ by $J$ from the right we do not change the product since $J^2=1$ in $PSL(2,\mathbb R)$ and transform the expression for $M$ into $M= H(x^{i_1})\mathbf E^\pm H(x^{i_2})\ldots H(x^{i_p})\mathbf E^\pm$, where $\mathbf E^+=\left(\begin{array}{cc}1&1\\0&1\end{array}\right)$,  $\mathbf E^-=\left(\begin{array}{cc}1&0\\1&1\end{array}\right)$ and $H(x) =\left(\begin{array}{cc}x^{1/2}&0\\0&x^{-1/2}\end{array}\right)$. Since all coefficients of these matrices are Laurent monomials with positive coefficients, the trace of any power of $M$ is a Laurent polynomial with positive coefficients.

Now consider the case of the space $\mathcal T^a(S)$. The reconstruction rule implies that the traces of the monodromy matrices are Laurent polynomials of the coordinates in integral powers of the coordinates $\{a_\alpha\}$. Positivity follows from the existence of the monomial map $p:\mathcal T^a(S)\rightarrow \mathcal T^x(S)$.

We leave the proof of the Laurent property and positivity for non closed curves as an exercise.

\appendix
\section{Combinatorial description of  $\mathcal D(S)$.} \label{combinatorial}\index{mapping class group}
 Denote by $|\Gamma|(S)$ the set of combinatorial types of triangulations of $S$. For each element of $|\Gamma|(S)$ fix a {\em marking}, i.e., a numeration of the edges of the triangulation. Denote by $\Gamma(S)$ the set of isotopy classes of marked triangulations of $S$. The presence of the marking changes the set of triangulations since some of them may have nontrivial symmetry group. Introducing the marking is a tool to remove this symmetry. (Here and below the vertical lines $|\cdot|$ indicate the diffeomorphism class.)

The mapping class group ${\cal D}(S)$ obviously acts freely on the space of marked triangulations having the space of combinatorial types of triangulations as a quotient,
$$
\mathbf\Gamma(S)/{\cal D}(S) = |\mathbf\Gamma|(S).
$$

Recall that a group can be thought of as a category with one object and invertible morphisms. Similarly, a {\em groupoid}\index{groupoid} is a category where all morphisms are isomorphisms, and any two objects are isomorphic. Since the automorphism groups of different objects of a groupoid are isomorphic, we can associate a group to a groupoid, well defined up to an isomorphism. We are going to construct a groupoid providing the mapping class group and admitting a simpler description by generators and relations than the mapping class group itself.

\begin{definition}
Let $|\mathbf \Gamma|(S)$ be the set of objects. For any two triangulations $|\Gamma|, |\Gamma_1| \in |\mathbf \Gamma(S)|$ let a morphism from $|\mathbf\Gamma|$ to $|\Gamma_1|$ be a pair of triangulations of $S$ of types $|\Gamma|$ and $|\Gamma_1|$ modulo the diagonal mapping class group action; we denote this morphism by $|\Gamma,\Gamma_1|$. For any three triangulations $\Gamma,\Gamma_1,\Gamma_2$, the composition of $|\Gamma,\Gamma_1|$ and $|\Gamma_1,\Gamma_2|$ is $|\Gamma,\Gamma_2|$. The described category is called the {\em modular groupoid}.
\end{definition}

One can easily verify that ({\rm 1}) the composition of morphisms is well defined; ({\rm 2}) the class of the pair of identical triangulations $|\Gamma,\Gamma|$ is the identity morphism and the inverse of the morphism $|\Gamma,\Gamma_1|$ is $|\Gamma_1,\Gamma|$; ({\rm 3}) the group of automorphisms of an object is the mapping class group ${\cal D}(S)$.

To give a description of the modular groupoid by generators and relations we need to introduce distinguished sets of morphisms called {\em flips} and {\em symmetries}. Recall that a morphism $|\Gamma,\Gamma_\alpha|$ is a flip if the triangulation $\Gamma_\alpha$ is obtained from $\Gamma$ by removing an edge $\alpha$ and replacing it by another diagonal of the arising quadrilateral. We use the notation $\Gamma_\alpha$  to emphasise the relation of this triangulation to the triangulation $\Gamma$. Note that for a given triangulation $\Gamma$, several marked embedded graphs may be denoted by $\Gamma_\alpha$ because no marking of $\Gamma_\alpha$ is indicated.

To each symmetry $\sigma$ of a triangulation $\Gamma$ we assign an automorphism $|\Gamma, \sigma \Gamma|$.

There is no canonical identification of edges of different triangulations even if a morphism between them is given.  However, two triangulations are related by a flip, we can introduce such an identification.  We exploit this identification and denote the corresponding edges of different graphs by the same letter if it is clear which sequence of flips relating these graphs is considered. To avoid confusion, note that this identification has nothing to do with the marking.

In this notation, the triangulation $\Gamma_{\alpha_1\cdots\alpha_n}$ is the graph obtained as a result of consecutive flips $\alpha_n, \ldots, \alpha_1$ of edges of a given graph $\Gamma$.

There are three kinds of relations between flips, which are satisfied for
any choice of marking of the triangulations entering the relations.

\begin{proposition}
A square of a flip is a symmetry\/{\rm:} if $|\Gamma_\alpha ,\Gamma|$ is a flip in an edge $\alpha$, then $|\Gamma, \Gamma_\alpha|$ is also a flip and\footnote{The notation {\bf R.n} indicates the number $n$ of graphs entering this relation.}

\noindent{\bf R.2.}\ \ \
$|\Gamma, \Gamma_\alpha| |\Gamma_\alpha, \Gamma|=1$.

Flips in disjoint edges commute\/{\rm:} if $\alpha$ and $\beta$ are two edges , then

\noindent{\bf R.4.} \ \ \
$|\Gamma_{\alpha\beta}, \Gamma_\alpha||\Gamma_\alpha, \Gamma| = |\Gamma_{\alpha\beta}, \Gamma_\beta||\Gamma_\beta, \Gamma|$.

Five consecutive flips in edges $\alpha$ and $\beta$ having one common vertex is the identity\/{\rm:} for such~$\alpha$ and~$\beta$, the triangulations $\Gamma_{\alpha\beta}$ and $\Gamma_{\beta\alpha}$ are related by a flip and

\noindent{\bf R.5.} $|\Gamma,\Gamma_\alpha||\Gamma_\alpha, \Gamma_{\beta\alpha}||\Gamma_{\beta\alpha},\Gamma_{\alpha\beta}|| \Gamma_{\alpha\beta},\Gamma_\beta||\Gamma_\beta,\Gamma| =1$.\index{pentagon relation}
\end{proposition}

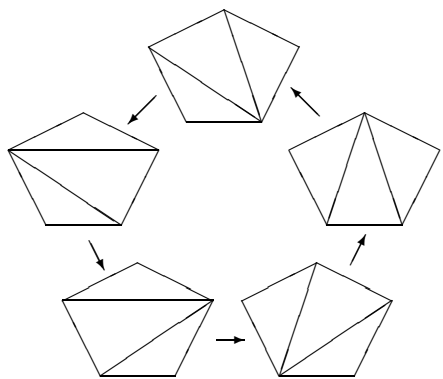
\begin{figure}
\unitlength 0.25mm
\begin{picture}(0,0)(-120,-145)
\put(0,0){\line(2,1){ 40}}
\put(40,20){\line(2,-1){ 40}}
\put(80,0){\line(-1,-2){ 20}}
\put(60,-40){\line(-1,0){40}}
\put(20,-40){\line(-1,2){ 20}}
\put(0,0){\line(1,0){80}}
\put(0,0){\line(3,-2){60}}
\put(43,-48){\vector(1,-2){8}}
\end{picture}
\begin{picture}(0,0)(-190,-200)
\put(0,0){\line(2,1){ 40}}
\put(40,20){\line(2,-1){ 40}}
\put(80,0){\line(-1,-2){ 20}}
\put(60,-40){\line(-1,0){40}}
\put(20,-40){\line(-1,2){ 20}}
\put(0,0){\line(3,-2){60}}
\put(40,20){\line(1,-3){20}}
\put(4,-26){\vector(-1,-1){15}}
\end{picture}
\begin{picture}(0,0)(-260,-145)
\put(0,0){\line(2,1){ 40}}
\put(40,20){\line(2,-1){ 40}}
\put(80,0){\line(-1,-2){ 20}}
\put(60,-40){\line(-1,0){40}}
\put(20,-40){\line(-1,2){ 20}}
\put(40,20){\line(1,-3){20}}
\put(40,20){\line(-1,-3){20}}
\put(16,18){\vector(-1,1){15}}
\end{picture}
\begin{picture}(0,0)(-230,-65)
\put(0,0){\line(2,1){ 40}}
\put(40,20){\line(2,-1){ 40}}
\put(80,0){\line(-1,-2){ 20}}
\put(60,-40){\line(-1,0){40}}
\put(20,-40){\line(-1,2){ 20}}
\put(40,20){\line(-1,-3){20}}
\put(20,-40){\line(3,2){60}}
\put(58,19){\vector(1,2){8}}
\end{picture}
\begin{picture}(0,230)(-130,-65)
\put(0,0){\line(2,1){ 40}}
\put(40,20){\line(2,-1){ 40}}
\put(80,0){\line(-1,-2){ 20}}
\put(60,-40){\line(-1,0){40}}
\put(20,-40){\line(-1,2){ 20}}
\put(20,-40){\line(3,2){60}}
\put(0,0){\line(1,0){ 80}}
\put(82,-21){\vector(1,0){15}}
\end{picture}
\caption{Pentagon relation.}
\label{pentagon}
\end{figure}

The proofs of the relations {\bf R.2} and {\bf R.4} are obvious. The relation {\bf R.5} is obvious from the triangulation of a pentagon shown on fig. \ref{pentagon}.

Replacing the mapping class group by the modular groupoid, we can express the latter through generators and relations in a simpler way.

Note that a symmetry can be represented as a ratio of two flips in a given edge and the modular groupoid is therefore generated by the flips only.
We do not describe relations between flips and graph symmetries in details because they are quite obvious. In fact, the symmetry groups of $\Gamma$ and $\Gamma_\alpha$ act transitively on the set of flips $|\Gamma, \Gamma_\alpha|$, and this action can be considered as relations between flips and symmetries.

Proposition~1 can be proved using direct combinatorial methods of the simplicial geometry. However, we give the main idea of another proof, which is more specific for the 2D situation.

\noindent{\em Proof of Proposition~1.}
To any connected simplicial complex ${\cal S}$ we can associate a groupoid by taking a point in each top-dimensional simplex for objects and  the homotopy classes of oriented paths connecting the chosen points as morphisms. The corresponding group is the fundamental group of the topological space given by the complex.

To any codimension one simplex we can associate two classes of paths (differing by orientation and having the identity morphism as their product) connecting adjacent top-dimensional simplices. It is natural to call them {\em flips}. We can associate a relation between the flips to any codimension two simplex. It is obvious that this set of flips generates the groupoid and that the only relation between the flips are given by codimension two simplices.

The same is true for an orbifold simplicial complex, where we replace simplices by quotients of simplices by finite groups. In this case, we must choose one generic point per each top-dimensional simplex as an object and orbifold homotopy classes of paths as morphisms. The corresponding group is the orbifold fundamental group of the orbifold given by the complex. The groupoid is now generated by flips and groups of top dimension simplices and still the only nontrivial relations are those given by codimension two simplices.

Consider now the Strebel \cite{Streb1} orbifold simplicial decomposition of the moduli space of complex structures on $S$. The orbifold fundamental group of the moduli space ${\cal M}$ is just the mapping class group ${\cal D}(S)$. Recall that Strebel orbisimplices are enumerated by ribbon graphs corresponding to~$S$ and the dimension of a simplex is equal to the number of its edges.  One can easily see that the groupoid of the Stre\-bel complex coincides with the modular groupoid. Moreover, the flips of the former correspond to the flips of the latter.  The relations between flips are given by codimension two cells, which correspond either to graphs with two four-valent vertices (which produces relation {\bf R.4}) or to graphs with one five-valent vertex (which produces relation {\bf R.5}). Relation {\bf R.2} holds true for any simplicial complex.\hfill{\bf q.e.d.}

\section{Markov numbers.}\index{Markov numbers}

 Consider a torus with one hole $T$.  The space of homotopy classes of simple (i.e., without self intersections) unoriented closed paths on it can be parameterised by points of ${\mathbb Q}P^1$. Indeed, once we have chosen an orientation of the path, we can consider it as an element of the first homology of $T$ with compact support. It is also obvious that any simple (indivisible) class is represented by a unique simple oriented closed path. Since the first homology group is ${\mathbb Z}^2$, it just gives the desired parameterisation.

Introduce the equiharmonic complex structure on $T$, i.e. the structure which has maximal symmetry group ${\mathbb Z}/3{\mathbb Z}$. For any closed path $\gamma$ on $T$ without self-intersections the numbers $X_\gamma = \frac{2}{3}\cosh l(\gamma)$, where $l(\gamma)$ are the lengths of the corresponding geodesic, are called {\em Markov numbers}.

The main properties of the Markov numbers are the following:

{\bf 1.} Markov numbers are positive integral.

{\bf 2.} Markov numbers include Fibonacci numbers with even indices
$2$, $5$, $13$, $34$, $89$, $233\ldots$.

Call a {\em Markov triple} a triple of Markov numbers $(\mathrm X,\mathrm Y,\mathrm Z)$ corresponding to three geodesic having pairwise one intersection point.

{\bf 3.} Elements of a Markov triple satisfy the Markov equation:
\begin{equation}
\mathrm X^2+\mathrm Y^2+\mathrm Z^2 = 3\mathrm X\mathrm Y\mathrm Z
\end{equation}

{\bf 4.} Any integer solution of this equation is a Markov triple.

{\bf 5.} For any Markov triple $(\mathrm X,\mathrm Y,\mathrm Z)$ the triples $(\mathrm Y,\mathrm Z,\mathrm X)$ and $(\mathrm Z,\mathrm Y-3\mathrm X\mathrm Z,\mathrm X)$ are also Markov triples. Any Markov triple can be obtained from the triple $(1,1,1)$ by a sequence of such transformations.

Since homotopy classes of closed non-selfintersecting curves can be parameterised by ${\mathbb Q}P^1$, one can choose an affine coordinate on ${\mathbb Q}P^1$ in such a way that the curves with coordinates $0,1$ and $\infty$ have
Markov numbers $1$. Denote by $M(u)$ the Markov number corresponding to the
curve with the coordinate $u \in {\mathbb Q}$.

{\bf 6.} The function $\psi(\frac{p}{q}) = \frac{1}{q}{\rm
arcosh}(\frac{3}{2}M(\frac{p}{q}))$, where $\gcd(p,q) = 1$, is extensible to
a continuous convex function on ${\mathbb R}$.

{\bf 7.} $M(x) = M(1-x) = M(\frac{1}{x}) = M(\frac{1}{1-x}) =
M(\frac{x}{x-1}) = M(\frac{x-1}{x})$

{\bf 8.} For any closed geodesics $\gamma$ on $S$ there exists a unique geodesics $\gamma\prime$ going from the puncture to the puncture which doesn't intersect $\gamma$. Let $l(\gamma\prime)$ be the length of the piece of $\gamma\prime$ between the intersection points with the horocycle surrounding the area $3$. Then $e^{l(\gamma\prime)} = M(\gamma).$

{\bf 9.} (Markov conjecture). The famous unproven Markov conjecture says that two Markov numbers $M(x)$ and $M(y)$ are different unless $x$ and $y$ are
related by transformations from property {\bf 7}. 

Taking into account that the segment $[0,1]$ is the fundamental domain of the action of transformations from property {\bf 7}, one can reformulate the Markov conjecture as that if $M(x) = M(y)$ and $x,y \in [0,1]$ then $x=y$.

{\em Proves of the properties.} (unfortunately, without the last one and the property {\bf 4}.)

There is only one combinatorial triangulation corresponding to the holed torus. It has one vertex, three edges and two triangles. This triangulation has obvious ${\mathbb Z}/3{\mathbb Z}$ symmetry group cyclically permuting the edges. Let $x,y,z$ be the corresponding coordinates on the Teichm\" uller space $\mathcal T^x(S)$.

 A closed curve on $S$ can be considered as a bounded lamination if we assign the weight $1$ to it. The standard coordinates of such laminations are given by three nonnegative integers $n_1,n_2,n_3$. These three numbers have no common factor, because otherwise the weight of the curve would be greater than $1$. On the other hand one of the numbers should be a sum of two others since otherwise there would be a component surrounding the hole. The relation between this parameterisation by $n_1,n_2,n_3$ and the parameterisation by ${\mathbb Q}P^1$ described above is given by
\begin{equation}
x =\left\{
\begin{array}{ll} \frac{-n_2}{n_1} & {\rm if }\  n_3 = n_1+n_2\\
                      \frac{n_2}{n_1} & {\rm if }\ n_1 = n_2+n_3\   {\rm or }\
                         n_2 = n_3+n_1
\end{array} \right.
\end{equation}

Denote by $\mathrm Z,\mathrm X$ and $\mathrm Y$ one thirds of traces of the elements of the Fuchsian group corresponding to the curves with co\-or\-di\-na\-tes $(1,1,0)$, $(0,1,1)$ and $(1,0,1)$ respectively. They can be easily computed using the explicit formulae for the Fuchsian group:
\begin{equation} \begin{array}{c}
\mathrm Z=\frac{1}{3}(x^{1/2}y^{1/2}+x^{1/2}y^{-1/2}+x^{-1/2}y^{-1/2}),\\
\mathrm X=\frac{1}{3}(y^{1/2}z^{1/2}+y^{1/2}z^{-1/2}+y^{-1/2}z^{-1/2}),\\ 
\mathrm Y=\frac{1}{3}(z^{1/2}x^{1/2}+z^{1/2}x^{-1/2}+z^{-1/2}x^{-1/2}).
\end{array}
 \label{mark1}
\end{equation}

Using these expressions we can verify the equality
\begin{equation}
\mathrm X^2 + \mathrm Y^2 + \mathrm Z^2 - 3\mathrm X \mathrm Y \mathrm Z =
-\frac{1}{9}(xyz - 2 + (xyz)^{-1})
\end{equation}

The symmetry of the graph obviously cyclically permutes the coordinates and
therefore the numbers $\mathrm Z,\mathrm X,\mathrm Y$. A flip of an edge acts by
the rule (\ref{al-flip}) and it results in the mapping

\begin{equation}
(\mathrm Z,\mathrm X,\mathrm Y) \mapsto
(\mathrm X,3\mathrm Y\mathrm Z-\mathrm X,\mathrm Z).\label{flipmark}
\end{equation}

If all three coordinates $x,y,z$ are ones, the corresponding complex
surface is just the equiharmonic punctured torus.

The properties 1,3,5,6 immediately follows from this picture. One can easily
check that $M(n)$ for $n \in {\mathbb N}$ are just the Fibonacci numbers what
gives the property 2. The property 7 is an immediate consequence of the
convexity property of the lamination length function. The property 4 was
proved by Markov himself.

The property 8 stands a little apart from the others since it is related to
the spaces ${\cal T}^a(S)$ and ${\mathsf T}^x(S)$ rather than ${\mathsf T}^a(S)$ and ${\cal T}^x(S)$ respectively. Consider a coordinate system $\mathrm U,\mathrm V,\mathrm W$ on ${\cal T}^a(S)$. Let  $\mathrm A= A_\rho$ be the area inside the only horocycle $\rho$. It easily follows from the expression for the area
\begin{equation}
 (\mathrm U^2 + \mathrm V^2 + \mathrm W^2) = \mathrm U\mathrm V\mathrm W\mathrm A  \label{markdual}
\end{equation}
The cyclic symmetry of the triangulation acts by cyclic permutation of $\mathrm U,\mathrm V,\mathrm W$. A flip of an edge acts by
\begin{equation}
(\mathrm U,\mathrm V,\mathrm W) \mapsto (\mathrm W,
\frac{\mathrm U^2+\mathrm W^2}{\mathrm V},\mathrm U).  \label{flipdecor}
\end{equation}

On the other hand this transformation law can be rewritten taking into account the equation (\ref{markdual}):
\begin{equation}
(\mathrm U,\mathrm V,\mathrm W) \mapsto (\mathrm W,\mathrm U\mathrm W\mathrm A - \mathrm V,\mathrm U)
\end{equation}

This rule coincides with (\ref{flipmark}) for $A=3$.

Now consider the decorated surface with $\mathrm U=\mathrm V=\mathrm W=1$. This is the surface with the area inside the horocycle $\mathrm A=3$. Applying modular transformations we get obviously the Markov triples, what proves the property 8.

There exists a canonical decomposition (called {\em Farey tessellation}) of the upper half plane $H$ into ideal triangles with vertices in all rational points of its ideal boundary. The dual graph to this tessellation is the universal three-valent tree. The faces of this tree are therefore in one-to-one correspondence with rational numbers. On the pictures below we have drawn a fragment of this tree with corresponding Markov numbers written on the faces.

\setlength{\unitlength}{0.175mm}
\begin{picture}(0,350)(70,500)
\put(400,780){\line( 0, 1){ 20}}
\put(400,800){\line( 1, 1){ 21}}
\put(420,820){\line( 1, 0){ 20}}
\put(400,800){\line(-1, 1){ 21}}
\put(420,820){\line(-1, 1){ 21}}
\put(240,700){\line( 4,-3){ 80}}
\put(320,640){\line( 1,-1){ 40}}
\put(360,600){\line( 1,-2){ 20}}
\put(160,640){\line( 1,-1){ 40}}
\put(200,600){\line( 1,-2){ 20}}
\put(320,640){\line(-1,-1){ 40}}
\put(280,600){\line(-1,-2){ 20}}
\put(280,600){\line( 1,-2){ 20}}
\put(360,600){\line(-1,-2){ 20}}
\put(200,600){\line(-1,-2){ 20}}
\put(120,600){\line( 1,-2){ 20}}
\put(400,780){\line(-2,-1){160}}
\put(240,700){\line(-4,-3){ 80}}
\put(160,640){\line(-1,-1){ 40}}
\put(120,600){\line(-1,-2){ 20}}
\put(400,780){\line( 2,-1){160}}
\put(560,700){\line( 4,-3){ 80}}
\put(640,640){\line( 1,-1){ 40}}
\put(680,600){\line( 1,-2){ 20}}
\put(680,600){\line(-1,-2){ 20}}
\put(640,640){\line(-1,-1){ 40}}
\put(600,600){\line(-1,-2){ 20}}
\put(600,600){\line( 1,-2){ 20}}
\put(560,700){\line(-4,-3){ 80}}
\put(480,640){\line(-1,-1){ 40}}
\put(440,600){\line(-1,-2){ 20}}
\put(440,600){\line( 1,-2){ 20}}
\put(480,640){\line( 1,-1){ 40}}
\put(520,600){\line( 1,-2){ 20}}
\put(520,600){\line(-1,-2){ 20}}
\put(100,560){\line(-2,-5){ 15.862}}
\put(100,560){\line( 2,-5){ 15.862}}
\put(140,560){\line(-2,-5){ 15.862}}
\put(140,560){\line( 2,-5){ 15.862}}
\put(180,560){\line(-2,-5){ 15.862}}
\put(180,560){\line( 2,-5){ 15.862}}
\put(220,560){\line(-2,-5){ 15.862}}
\put(220,560){\line( 2,-5){ 15.862}}
\put(260,560){\line(-2,-5){ 15.862}}
\put(260,560){\line( 2,-5){ 15.862}}
\put(300,560){\line(-2,-5){ 15.862}}
\put(300,560){\line( 2,-5){ 15.862}}
\put(340,560){\line(-2,-5){ 15.862}}
\put(340,560){\line( 2,-5){ 15.862}}
\put(380,560){\line(-2,-5){ 15.862}}
\put(380,560){\line( 2,-5){ 15.862}}
\put(420,560){\line(-2,-5){ 15.862}}
\put(420,560){\line( 2,-5){ 15.862}}
\put(460,560){\line(-2,-5){ 15.862}}
\put(460,560){\line( 2,-5){ 15.862}}
\put(500,560){\line(-2,-5){ 15.862}}
\put(500,560){\line( 2,-5){ 15.862}}
\put(540,560){\line(-2,-5){ 15.862}}
\put(540,560){\line( 2,-5){ 15.862}}
\put(580,560){\line(-2,-5){ 15.862}}
\put(580,560){\line( 2,-5){ 15.862}}
\put(620,560){\line(-2,-5){ 15.862}}
\put(620,560){\line( 2,-5){ 15.862}}
\put(660,560){\line(-2,-5){ 15.862}}
\put(660,560){\line( 2,-5){ 15.862}}
\put(700,560){\line(-2,-5){ 15.862}}
\put(700,560){\line( 2,-5){ 15.862}}
\put(400,820){\makebox(0,0)[cb]{1}}
\put(400,680){\makebox(0,0)[cb]{5}}
\put(240,640){\makebox(0,0)[cb]{29}}
\put(560,640){\makebox(0,0)[cb]{13}}
\put(160,600){\makebox(0,0)[cb]{169}}
\put(320,600){\makebox(0,0)[cb]{433}}
\put(480,600){\makebox(0,0)[cb]{194}}
\put(640,600){\makebox(0,0)[cb]{34}}
\tiny
\put(120,560){\makebox(0,0)[cb]{985}}
\put(200,560){\makebox(0,0)[cb]{14701}}
\put(280,560){\makebox(0,0)[cb]{37666}}
\put(360,560){\makebox(0,0)[cb]{6466}}
\put(440,560){\makebox(0,0)[cb]{2897}}
\put(520,560){\makebox(0,0)[cb]{7561}}
\put(600,560){\makebox(0,0)[cb]{1325}}
\put(680,560){\makebox(0,0)[cb]{89}}
\normalsize
\put(140,760){\makebox(0,0)[cb]{2}}
\put(660,760){\makebox(0,0)[cb]{1}}
\tiny
\put(100,510){\makebox(0,0)[cb]{5741}}
\put(140,500){\makebox(0,0)[cb]{499393}}
\put(180,510){\makebox(0,0)[cb]{7453378}}
\put(220,500){\makebox(0,0)[cb]{1278818}}
\put(260,510){\makebox(0,0)[cb]{3276569}}
\put(300,500){\makebox(0,0)[cb]{48928105}}
\put(345,510){\makebox(0,0)[cb]{8399329}}
\put(380,500){\makebox(0,0)[cb]{96557}}
\put(420,510){\makebox(0,0)[cb]{43261}}
\put(460,500){\makebox(0,0)[cb]{1686049}}
\put(500,510){\makebox(0,0)[cb]{4400489}}
\put(540,500){\makebox(0,0)[cb]{294685}}
\put(580,510){\makebox(0,0)[cb]{51641}}
\put(620,500){\makebox(0,0)[cb]{135137}}
\put(660,510){\makebox(0,0)[cb]{9077}}
\put(700,500){\makebox(0,0)[cb]{233}}
\end{picture}

As a concluding remark of this section note that, as it was observed by
A.Bondal, Markov triples are dimensions of elements of distinguished sets of
sheaves on ${\mathbb C}P^2$. The relations between these two ways of obtaining
Markov numbers are completely unclear and very exciting.

\vspace{2mm}
\noindent V.V.F: Institute for Theoretical and Experimental Physics, \\ B.Cheremushkinskaya 25, 117259 Moscow, Russia {\tt  fock@math.brown.edu}

\noindent A.B.G: Brown University, \\ Box 1917, 02906, Providence, RI, USA. {\tt  sasha@math.brown.edu}

\end{document}